\documentclass[12pt]{article}
\usepackage{amssymb}
\usepackage{extarrows}
\usepackage{amsmath}
\usepackage{cite}
\newtheorem{theorem}{Theorem}[section]
\newtheorem{definition}[theorem]{Definition}
\newtheorem{lemma}[theorem]{Lemma}
\newtheorem{corollary}[theorem]{Corollary}
\newtheorem{proposition}[theorem]{Proposition}
\newtheorem{example}[theorem]{Example}
\newtheorem{remark}[theorem]{Remark}
\newtheorem{theorem(Composition-Diamond Lemma)}[theorem]{Theorem(Composition-Diamond Lemma)}

\begin{document}

\title{A new Composition-Diamond lemma for associative conformal algebras\footnote{Supported by the NNSF of China (11171118, 11571121).}}

\author{
Lili Ni and Yuqun
Chen\footnote {Corresponding author.}  \\
{\small \ School of Mathematical Sciences, South China Normal
University}\\
{\small Guangzhou 510631, P. R. China}\\
{\small nilili2009@163.com}\\
{\small yqchen@scnu.edu.cn}}

\date{}
\maketitle

\maketitle \noindent\textbf{Abstract:}
Let $C(B,N)$ be the free associative conformal algebra
generated by a set $B$ with a bounded locality $N$.
Let $S$ be a subset of $C(B,N)$. A Composition-Diamond lemma for
associative conformal algebras is firstly established
by  Bokut, Fong, and Ke in 2004  \cite{BFK04} which
claims that if (i) $S$ is a Gr\"{o}bner-Shirshov basis
in $C(B,N)$, then (ii) the set of $S$-irreducible words
is a linear basis of the quotient conformal algebra
$C(B,N|S)$, but not conversely. In this paper,
by introducing some new definitions of normal $S$-words,
compositions and compositions to be trivial, we give a
new Composition-Diamond lemma for associative conformal
algebras which makes the conditions  (i) and (ii)
 equivalent. We show that for each ideal $I$ of $C(B,N)$, $I$ has
  a unique reduced Gr\"{o}bner-Shirshov basis. As applications, we show that
 Loop Virasoro Lie conformal algebra and Loop Heisenberg-Virasoro
 Lie conformal algebra are embeddable into their universal
enveloping associative conformal algebras.

\ \


\ \

\noindent {\bf AMS} Mathematics Subject Classification(2000): 17B69,
16S15, 13P10,  08A50

\section{Introduction}

The subject of conformal algebras is closely related to vertex
algebras (see,  V. Kac \cite{Kac96}). Implicitly, vertex algebras
were introduced by Belavin, Polyakov, and Zamolodchikov in 1984
\cite{Belavin}. Explicitly, the definition of vertex algebras was
given by R. Borcherds in 1986 \cite{Bor86}, which led to his
solution of the Conway-Norton conjecture in the theory of finite
simple groups \cite{Bor92, Frenkel88}.  As pointed out by Kac
\cite{Kac97, Kac96}, conformal and vertex algebras provide a
rigorous mathematical study of the ¡°locality axiom¡± which came
from Wightman¡¯s axioms of quantum field theory \cite{Wightman}. M.
Roitman studied free (Lie and associative) conformal and vertex
algebras in \cite{Ro99}. Free vertex algebras were mentioned in the
original paper of Borcherds \cite{Bor86}. Since conformal and vertex
algebras are not varieties in the sense of universal algebra (see,
P.M. Cohn \cite{Cohn81}), the existence of free conformal and free
vertex algebras is not guaranteed by the general theory and should
be proved. It was done by Roitman  \cite{Ro99}. The free associative
conformal algebra generated by a set $B$ with a  locality function
$N(-,-):\ B\times B\rightarrow \mathbb{Z}_{\geq 0}$ is constructed
by Bokut, Fong, and Ke in 1997 \cite{BFK0}.

 Gr\"{o}bner bases and Gr\"{o}bner-Shirshov bases  were invented
independently by A.I. Shirshov for ideals of free (commutative,
anti-commutative) non-associative algebras {\cite{Shir3,Sh62a}}, free
Lie algebras \cite{Shir3} and implicitly free associative
algebras \cite{Shir3}  (see also \cite{Be78,Bo76}), by H.
Hironaka \cite{Hi64} for ideals of the power series algebras (both
formal and convergent), and by B. Buchberger \cite{Bu70} for ideals
of the polynomial algebras.
Gr\"{o}bner bases and Gr\"{o}bner-Shirshov bases theories have been
proved to be very useful in different branches of mathematics,
including commutative algebra and combinatorial algebra. It is a
powerful tool to solve the following classical problems: normal
form; word problem; conjugacy problem; rewriting system; automaton;
embedding theorem; PBW theorem;  extension; homology;  growth
function; Dehn function; complexity; etc. See, for example, the
books \cite{AL, BKu94, BuCL, BuW, CLO, Ei} and the surveys \cite{BC,BC13, BFKK00, BK03, BK05}.

In Gr\"{o}bner-Shirshov bases theory for a category of algebras,
the key part is to establish ``Composition-Diamond lemma" for such algebras. The name
``Composition-Diamond lemma"
combines the Neuman Diamond Lemma \cite{Newman},
the Shirshov Composition Lemma \cite{Sh62a}
and the Bergman Diamond Lemma \cite{Be78}.

Gr\"{o}bner-Shirshov bases theory for associative conformal
algebras is firstly established by  Bokut, Fong, and Ke in 2004  \cite{BFK04},
i.e., they give a Composition-Diamond lemma for associative conformal algebras.

Let $C(B,N)$ be the free associative conformal algebra over a field $\mathbf{k}$ of characteristic 0
generated by a set $B$ with a bounded locality $N$. Let $S\subset C(B,N)$
be a monic subset of polynomials and $Id(S)$ be
the ideal of $C(B,N)$ generated by $S$. A normal word $[u]$ is said to
be $S$-irreducible if $u\neq \overline{[v]_{_{D^is}}}$, where $[v]_{_{D^is}}$
is any normal $S$-word. Let $Irr(S)$ be the set of all $S$-irreducible words.
Consider the following statements:
\begin{enumerate}
 \item[(i)]\ The set $S$ is a Gr\"{o}bner-Shirshov basis in $C(B,N)$.

 \item[(ii)]\ $Irr(S)$ is a $\mathbf{k}$-basis of $C(B,N|S)=C(B,N)/Id(S)$.
\end{enumerate}

In \cite{BFK04}, it is shown that $(i)\Rightarrow(ii)$ but $(ii)\not\Rightarrow(i)$.

In this paper, by introducing some new definitions of normal $S$-words,
compositions and compositions to be trivial,
we give a new Composition-Diamond lemma for associative conformal algebras
which makes the two conditions above
 equivalent, see Theorem \ref{t1}. We show that for each ideal $I$ of $C(B,N)$,
 $I$ has a unique reduced Gr\"{o}bner-Shirshov basis. Of course, our definition of a
 Gr\"{o}bner-Shirshov basis in $C(B,N)$ is different from one in  \cite{BFK04}.
As applications, we give a characterization of a conformal
 algebra $C(B,N|S)$ to have a $\mathbf{k}[D]$-basis under some
 conditions, see Proposition \ref{t2}, and PBW theorems:
 Loop Virasoro Lie conformal algebra and Loop Heisenberg-Virasoro
 Lie conformal algebra  are embeddable into their universal
enveloping associative conformal algebras, respectively,
see Examples \ref{exam1} and \ref{exam2}.

We are grateful to Professor  L.A. Bokut for  guidance
and to Mr. Zerui Zhang for useful discussions.

\section{Preliminary}

We begin with the formal definition of a conformal algebra.

\begin{definition}(\cite{BFK0,BFKK00,BFK00})
A conformal algebra $C=(C,(n),n\in\mathbb{Z}_{\geq 0},D)$ is a
linear space over a field $\mathbf{k}$ of characteristic $0$,
equipped with bilinear multiplications $a_{(n)}b$,
$n\in\mathbb{Z}_{\geq 0}=\{0,1,2,\cdots \}$,
and a linear map $D$, such that the following axioms are valid:
\begin{enumerate}
\item[(C1)]\ (locality) For any $a, b\in C$, there exists
a nonnegative integer $N(a,b)$ such that $a_{(n)}b=0$ for
$n\geq N(a,b)$\ ($N(a, b)$ is called the ordering of locality of $a$ and $b$).
\item[(C2)]\ $D(a_{(n)}b)=Da_{(n)}b+a_{(n)}Db$ for any $a,b\in C$
and $n\in\mathbb{Z}_{\geq 0}$.
\item[(C3)]\ $Da_{(n)}b=-na_{(n-1)}b$ for any $a,b\in C,\ n>0$ and $Da_{(0)}b=0$.
\end{enumerate}

A conformal algebra $C$ is called  {\it associative}
if the following identity holds for
all $a,b,c\in C,\ m,n\in\mathbb{Z}_{\geq 0}$,

$\bullet$ (Associative identity)
$$
(a_{(n)}b)_{(m)}c=\sum_{t\geq 0}(-1)^{t}\binom{n}{t}a_{(n-t)}(b_{(m+t)}c).
$$

A conformal algebra $L=\langle L, [n], n\in\mathbb{Z}_{\geq 0},
D\rangle$  is called a Lie conformal algebra if $L$ satisfies the
following two axioms:

$\bullet$ (Anti-commutativity)\ $a_{[n]}b=-\{b_{[n]}a\}$, where
$$
\{b_{[n]}a\}=\sum\limits_{k\geq
0}(-1)^{n+k}\frac{1}{k!}D^{k}(b_{[n+k]}a).
$$

$\bullet$ (Jacobi identity)
$$
(a_{[n]}b)_{[m]}c=\sum\limits_{k\geq 0}(-1)^{k}
\binom{n}{k}(a_{[n-k]}(b_{[m+k]}c)-b_{[m+k]}(a_{[n-k]}c)).
$$
\end{definition}

Let $B$ be a set. An associative conformal algebra
$$
C(B,N)=( B;\ a_{(n)}b=0,\ a,b\in B,\ n\geq N(a,b))
$$
is called the free associative conformal algebra
generated by $B$ with the locality function $N:B\times
B\longrightarrow \mathbb{Z}_{\geq 0}$,  if for
any associative conformal algebra $C$ and any maping $\varepsilon :
B\longrightarrow C$ with
$\varepsilon (a)_{(n)}\varepsilon (b)=0$ for all $a,b\in B$ and
$n\geq N(a,b)$, there exists a unique  homomorphism
$\varphi: C(B,N)\rightarrow C$ such that the following diagram
is commutative:

\setlength{\unitlength}{1cm}
\begin{picture}(7, 3)
\put(4.2,2.3){\vector(1,0){1.7}} \put(4.1, 2.0){\vector(0,-1){1.3}}
\put(5.9,2.1){\vector(-1,-1){1.6}} \put(3.9,0.2){$C$}
\put(6,2.2){$C(B,N)$} \put(3.9,2.2){$B$} \put(4.9, 2.4){$i$}
\put(3.7,1.3){$\varepsilon$} \put(5.3,1){$\exists!\varphi$}
\end{picture}

\noindent where $i$ is the inclusion map.

Let $D^{\omega}(B)=\{D^{i}b\mid b\in B,\ i\geq 0\}$.
We define (non-associative) words on $D^{\omega}(B)$ inductively:
\begin{enumerate}
\item[i)]\ For any $b\in B$ and $i\geq 0$,  $D^ib$ is
a word of length 1.
\item[ii)]\ If $(u)$ and $(v)$ are two words of lengths
$k$ and $l$, respectively,
then $(u)_{(n)}(v)$ is a word of length $k+l$ for any
$n\in \mathbb{Z}_{\geq 0}$.
\end{enumerate}
For any word $(u)$ on $D^{\omega}(B)$, we denote its length by $|u|$.
A word $(u)$ is a normal word if it has the form:
$$
(u)=b_{_{1}(n_1)}(b_{_{2}(n_2)}\cdots (b_{_{k}(n_k)}D^{i}b_{_{k+1}})\cdots )\ \ \ \ \mbox{(right normed bracketing)},
$$
where $b_{_{k+1}}, b_{_j}\in B,\ 0\leq n_{_j}< N(b_{_j}, b_{_{j+1}}),\
1\leq j\leq k,\ k\geq 0,\ i\geq 0$. If this is the case, then
denote by $ind(u)=i$.

\begin{lemma} (\cite{BFK0}, Lemma 2.3)\label{l0}
Any word $(u)$ (on $D^{\omega}(B)$) is a linear combination
of normal words of the same length $|u|$.
\end{lemma}

The free associative conformal algebra $C(B,N)$
generated by a set $B$ with a locality function $N(-,-)$
is constructed in the paper \cite{BFK0}.
In \cite{BFK0,Ro99}, it was proved that there is
a linear basis of $C(B,N)$ consisting of
normal words, see \cite{BFK0}, Theorem 2.5.

Gr\"{o}bner-Shirshov bases theory for free associative conformal algebras
 with a uniform bounded locality function $N$ is established in the paper
 \cite{BFK04}. In this paper, we only consider the situation when
 the locality function $N(b,b'), b,b'\in B$ is uniformly bounded
 by some number $N$. Then without loss of generality,
 we assume that $N(b,b')=N$ for all $b,b'\in B$.

Let $C(B,N)$ be the free associative conformal algebra
generated by a set $B$ with a bounded locality $N$.
A normal word is an expression of the form
\begin{equation}\label{e1.1}
[u]=b_{_{1}(n_1)}(b_{_{2}(n_2)}\cdots (b_{_{k}(n_k)}D^{i}b_{_{k+1}})\cdots)
\end{equation}
where $b_{k +1},b_{_j}\in B,\ 0\leq n_{_j}< N,\
1\leq j\leq k,\ k,i\geq0$ and $|u|=k+1$.
Here and after, $[u]$ would mean the right normed bracketing.
A normal word without brackets
\begin{equation}\label{e1.2}
u=b_{_{1}(n_1)}b_{_{2}(n_2)}\cdots b_{_{k}(n_k)}D^{i}b_{_{k+1}}
\end{equation}
is referred to as an associative  normal word.
Denote
$$
T:=\{b_{_{1}(n_1)}b_{_{2}(n_2)}\cdots b_{_{k}(n_k)}D^{i}b_{_{k+1}}\mid
b_{_{k+1}}, b_{_j}\in B,\ 0\leq n_{_j}< N,\
1\leq j\leq k,\ k,i\geq 0\}.
$$
Let $u=b_{_{1}(n_1)}b_{_{2}(n_2)}\cdots b_{_{k}(n_k)}D^{i}b_{_{k+1}},\
v=b'_{_{1}(n'_1)}b'_{_{2}(n'_2)}\cdots b'_{_{t}(n'_t)}D^{j}b'_{_{t+1}}\in T$.
Denote
\begin{eqnarray*}
uD^l&:=&b_{_{1}(n_1)}b_{_{2}(n_2)}\cdots b_{_{k}(n_k)}D^{i+l}b_{_{k+1}},\\
u^{\backslash D}&:=&b_{_{1}(n_1)}b_{_{2}(n_2)}\cdots b_{_{k}(n_k)}b_{_{k+1}},\\
u^{\backslash D}\natural v&:=&b_{_{1}(n_1)}b_{_{2}(n_2)}\cdots b_{_{k}(n_k)}
b_{_{k+1}(N-1)}b'_{_{1}(N-1)}b'_{_{2}(N-1)}\cdots b'_{_{t}(N-1)}D^{j}b'_{_{t+1}}.
\end{eqnarray*}

Any polynomial $f\in C(B,N)$ is
a linear combination of normal words.
We shall refer to normal word $[u]$ as D-free if $ind(u)=0$,
and say that $f$ is D-free if every monomial word in $f$ is D-free.
Let $S\subseteq C(B,N)$.
Then $S$ is said to be $D$-free if all polynomials in $S$ are $D$-free.

\section{Composition-Diamond lemma for associative conformal algebras}

Let $B$ be a well-ordered set.
We order the words of the forms (\ref{e1.1}) and (\ref{e1.2}) according to
the lexicographical ordering of their weights (see \cite{BFK04}).

Let $u=b_{_{1}(n_1)}b_{_{2}(n_2)}\cdots b_{_{k}(n_k)}D^{i}b_{_{k+1}}\in T$. Denote
$$
wt(u)=(|u|,b_1,n_1,\cdots ,b_k,n_k, b_{k+1},i).
$$
Then for any $u, v\in T$, $[u]>[v]\Leftrightarrow u>v$, where
\begin{equation}\label{e1.3}
u>v\Leftrightarrow wt(u)>wt(v)\ \ \ \ \mbox{ lexicographically.}
\end{equation}

It is clear that such an ordering is a well ordering on $T$.
We will use this ordering in the sequel.
For convenience, we define $u>0$ for any $u\in T$.

Given $f\in C(B,N)$,  we denote by $\bar{f}$ the leading
associative normal word in $f$ with respect to the ordering (\ref{e1.3})
and say $f$  monic if the coefficient of $\bar{f}$ is 1.
In particular,  we denote $\bar{0}=0$. A subset $S$ of $C(B,N)$ is said monic if
for any $s\in S$, $s$ is monic.

\begin{lemma}(\cite{BFK04})\label{1}
Suppose that $a\in T$ is $D$-free, $0\leq n<N$ and $u, v\in T$. Then
\begin{enumerate}
 \item[(i)]\ $\overline{[a]_{(n)}[u]}=a_{(n)}u$.
 \item[(ii)]\ $u>v\Longrightarrow \overline{[a]_{(n)}[u]}>\overline{[a]_{(n)}[v]}$.
 \item[(iii)]\ $\overline{D^{i}([u])}=uD^{i}$, $i\in \mathbb{Z}_{\geq 0}$.
 \item[(iv)]\ $u>v\Longrightarrow \overline{D^{i}([u])}>
 \overline{D^{i}([v])}$, $ i\in \mathbb{Z}_{\geq 0}$.
\end{enumerate}
\end{lemma}

Note that the ordering $>$ on $T$ is not compatible with the multiplications,
for example, in $C(\{a\},N=1)$, we have
$a_{(0)}Da>a_{(0)}a$, but
$$
0=(a_{(0)}Da)_{(0)}a< (a_{(0)}a)_{(0)}a=a_{(0)}(a_{(0)}a).
$$

\begin{lemma}\label{l2}
Let $[u]$ and $[v]$ be normal words. Then for any $n\in \mathbb{Z}_{\geq 0}$,
$\overline{[u]_{(n)}[v]}\leq u^{^{\backslash D}} \natural v$.
\end{lemma}
\noindent{\bf Proof.}  Induction on $(|u|, |v|, n)$, where the ordering is lexicographically.

Case 1\ \  $|u|=|v|=1$. Suppose $[u]=D^ib$ for some $b\in B$ and $i\geq 0$.
Then $D^ib_{(n)}[v]=(-1)^{i}\frac{n!}{(n-i)!}b_{(n-i)}[v]$.
So it suffices to consider $i=0$.
If $0\leq n<N$, then $b_{(n)}[v]$ is a normal word.
Thus we have $b_{(n)}[v]\leq b\natural v$.
Now let $n\geq N$. Assume that $v=D^jb'$ for some $b'\in B,\ j\geq 0$.
If $j=0$, then $b_{(n)}b'=0$. If $j>0$, then
$
b_{(n)}D^jb'=-\sum_{k\geq 1}(-1)^{k}
\binom{j}{k}\frac{n!}{(n-k)!}b_{(n-k)}D^{j-k}b'.
$
We can get $\overline{b_{(n)}D^jb'}\leq b \natural D^jb'$
by induction on $n$.

This shows that the result holds for $(|u|, |v|, n)=(1,1,n)$ for any $n\in \mathbb{Z}_{\geq 0}$.

Case 2\ \  $|u|=1$ and $|v|>1$. Assume that $[v]=b'_{(m)}[v_{_1}]$.
Then
$$
b_{(n)}(b'_{(m)}[v_{_1}])=
-\sum\limits_{k\geq 1}(-1)^{k}\binom{n}{k}b_{(n-k)}(b'_{(m+k)}[v_{_1}]).
$$
By induction, we have
$\overline{b'_{(t)}[v_{_1}]}\leq b' \natural v_{_1} $ for any $t\geq 0$
and
$\overline{b_{(n)}(b'_{(m)}[v_{_1}])}\leq b_{(N-1)}b' \natural v_{_1}=b \natural v$.

Case 3\ \  $|u|>1$. Let $[u]=b_{(m)}[u_{_1}]$ for some
$b\in B, u_{_1}\in T$. Then
$$
[b_{(m)}u_{_1}]_{(n)}[v]=\sum\limits_{k\geq 0}(-1)^{k}
\binom{m}{k}b_{(m-k)}([u_{_1}]_{(n+k)}[v]).
$$
By induction, we have $\overline{[u_{_1}]_{(n+k)}[v] }\leq
u_{_1}^{^{\backslash D}} \natural v$.
Hence
$
\overline{[u]_{(n)}[v]}= b_{(m)}(\overline{[u_{_1}]_{(n)}[v]})
\leq b_{(m)}u_{_1}^{^{\backslash D}} \natural v
=u^{^{\backslash D}} \natural v.
$

This completes the proof.
\ \ \ \  $\square$

\begin{corollary}\label{c1}
Suppose that $u, v\in T,\ ind(u)=0,\ u>v$ and $0\leq n<N$.
Then for any $w\in T$, the following statements hold.
\begin{enumerate}
 \item[(i)]\ $u_{(n)}w>v^{\backslash D}\natural  w$.
 \item[(ii)]\ $\overline{[u]_{(n)}[w]}>\overline{[v]_{(n)}[w]}$.
 \end{enumerate}
\end{corollary}
\noindent{\bf Proof.}
(i)\ Since $u>v$ and $u$ is $D$-free, we have $u>v^{^{\backslash D}}$.
Thus $u_{(n)}w>v^{^{\backslash D}}\natural w$ by the definition of
ordering (\ref{e1.3}).

(ii)\ By Lemma \ref{l2}, we have
$\overline{[v]_{(n)}[w]}\leq v^{^{\backslash D}}\natural w$.
By Lemma \ref{1}, we have $\overline{[u]_{(n)}[w]}=u_{(n)}w$.
Hence $\overline{[u]_{(n)}[w]}>\overline{[v]_{(n)}[w]}$.
\ \ \ \  $\square$

\begin{corollary}\label{c2}
For any $n\geq 0,\ u\in T,\ f\in C(B,N)$, we have
$\overline{f_{(n)}[u]}\leq \bar{f}^{^{\backslash D}}\natural u$.
\end{corollary}
\noindent{\bf Proof.} If $f=0$ then the result is clear.
Let $0\neq f=\alpha[\bar{f}]+\sum\alpha_{_t}[v_{_t}]$,
where $\alpha,\alpha_{_t}\in \mathbf{k},\ \bar{f}, v_{_t}\in T$
and $v_{_t}<\bar{f}$. Then we have
$
f_{(n)}[u]=\alpha[\bar{f}]_{(n)}[u]+\sum\alpha_{_t}[v_{_t}]_{(n)}[u]
$
and $ v_{t}^{^{\backslash D}}\leq\bar{f}^{^{\backslash D}}$.
Thus $v_{t}^{^{\backslash D}}\natural u\leq\bar{f}^{^{\backslash D}}\natural u$.
By Lemma \ref{l2}, we have
$\overline{[\bar{f}]_{(n)}[u]}\leq \bar{f}^{^{\backslash D}}\natural u$
and $\overline{[v_{_t}]_{(n)}[u]}\leq v_{t}^{^{\backslash D}}\natural u$.
Hence $\overline{f_{(n)}[u]}\leq \bar{f}^{^{\backslash D}}\natural u$.
\ \ \ \  $\square$

\ \

Now we define $S$-words and normal $S$-words.

Let $S\subset C(B,N)$ be a set of monic polynomials. We define $S$-words
$(u)_{_{D^is}}$ by induction.

(i)\ $(D^is)_{_{D^is}}=D^is$, where $s\in S$ and $i\geq 0$,
is an $S$-word of $S$-length 1.

(ii)\ If $(u)_{_{D^is}}$ is an $S$-word of $S$-length $k$,
and $(v)$ is any word  of length $l$,
then $(u)_{_{D^is}(m)}(v)$ and $(v)_{(m)}(u)_{_{D^is}}$
are $S$-words of $S$-length $k+l$.

The $S$-length of an $S$-word $(u)_{_{D^is}}$ will be denoted by $|u|_{_{D^is}}$.

\begin{definition}
Let $S\subset C(B,N)$ be a set of monic polynomials.
An associative normal $S$-word is an expression of the following forms:
\begin{eqnarray}\label{e1.4}
 u_s=a_{(n)}s_{(m)}c
\end{eqnarray}
where $s\in S$, $a,c\in T$ ($a$ may be empty) with $a$ and $\bar{s}$ being D-free,
$0\leq n,m<N$; or
\begin{eqnarray}\label{e1.5}
u_{_{D^is}}=a_{(n)}D^is
\end{eqnarray}
where $s\in S, i\geq0, a\in T$ ($a$ may be empty)  is D-free, $0\leq n<N$.
\end{definition}
Note that in (\ref{e1.4}) and (\ref{e1.5}), if $a$ is empty, then (\ref{e1.4}) means
$u_s=s_{(m)}c$ and (\ref{e1.5}) means $u_{_{D^is}}=D^is$.
From now on, $ a=1$ means that $a$ is empty.

For convenience, we call associative normal $S$-words (\ref{e1.4}) and (\ref{e1.5}) the first kind
and the second kind associative normal $S$-words, respectively, while
\begin{eqnarray}\label{e1.5.1}
[u]_{_s}=[a_{(n)}s_{(m)}c] \ \ \mbox{ and } \ \
[u]_{_{D^{i}s}}=[a_{(n)}D^is]
\end{eqnarray}
will be referred to normal $S$-words of the first kind
and the second kind, respectively.
A common notation for  (\ref{e1.5.1}) is
$[u]_{_{D^is}}$ as both two kinds of normal $S$-words, i.e.,
\begin{eqnarray*}
[u]_{_{D^is}}=\left\{
\begin{aligned}
& [s_{(m)}c],\ \ \ \ \ \ \ \ \ \ \ ind(\bar{s})=0,   \\
&[a_{(n)}s_{(m)}c],\ \ \ \ \ \    ind(\bar{s})=0, \\
&D^is,\ \ \ \ \ \     \\
&[a_{(n)}D^is],\ \ \ \ \ \
\end{aligned}
\right.
\end{eqnarray*}
where $s\in S$, $a,c\in T$, $a$ is D-free,
$ i\geq 0,\ 0\leq n,m<N$.

From now on, we always assume that $S\subset C(B,N)$ is a set of monic polynomials.

\begin{lemma}\label{l3}
Let $[u]_{_{D^is}}$ be a normal $S$-word.
Then the following statements hold.
\begin{enumerate}
 \item[(i)]\ If $[u]_{_{D^is}}=[a_{(n)}s_{(m)}c]$, then
 $\overline{[a_{(n)}s_{(m)}c]}=a_{(n)}\bar{s}_{(m)}c$.
 \item[(ii)]\ If $[u]_{_{D^is}}=[a_{(n)}D^is]$,
 then $\overline{[a_{(n)}D^is]}=a_{(n)}\bar{s}D^i$ .
\end{enumerate}
\end{lemma}
\noindent{\bf Proof.}
Let $s=[\bar{s}]+\sum\alpha_{_t}[u_{_t}]$, where
$\bar{s}, u_{_t}\in T$
and $u_{_t}<\bar{s}$. Then

(i)\ $[a_{(n)}s_{(m)}c]=[a_{(n)}[\bar{s}]_{(m)}c]+
\sum\alpha_{_t}[a_{(n)}[u_{_t}]_{(m)}c]$.
By Corollary \ref{c1}, we have
$\overline{[u_{_t}]_{(m)}[c]}<\overline{[\bar{s}]_{(m)}[c]}$.
By Lemma \ref{1}, $\overline{[\bar{s}]_{(m)}[c]}=\bar{s}_{(m)}c$.
Therefore $\overline{[a_{(n)}s_{(m)}c]}=a_{(n)}\bar{s}_{(m)}c$.

(ii)\ $[a_{(n)}D^is]=[a_{(n)}D^{i}[\bar{s}]]+
\sum\alpha_{_t}[a_{(n)}D^{i}[u_{_t}]]$.
By Lemma \ref{1}, we have $\overline{D^{i}[\bar{s}]}=\bar{s}D^i$ and
$\overline{D^{i}[u_{_t}]}=u_{_t}D^i$.
Since $u_{_t}<\bar{s}$, we have $u_{_t}D^i<\bar{s}D^i$.
Hence $\overline{[a_{(n)}D^is]}=a_{(n)}\bar{s}D^i$.
\ \ \ \  $\square$

\begin{lemma}\label{l4}
Let $[u]_{_{D^is}}$ be a normal $S$-word.
Then for any $j>0$, the following statements hold.
\begin{enumerate}
 \item[(i)]\ If $[u]_{_{D^is}}=[a_{(n)}s_{(m)}c]$, then
 $D^j[a_{(n)}s_{(m)}c]=[a_{(n)}s_{(m)}c D^j]
 +\sum\alpha_{_t}[a_{_{t}(n_t)}s_{(m_t)}c_t]$, where
each $[a_{_{t}(n_t)}s_{(m_t)}c_t]$ is a normal $S$-word,
 $a_{_{t}(n_t)}\bar{s}_{(m_t)}c_t<a_{(n)}\bar{s}_{(m)}c D^j$.
 \item[(ii)]\ If $[u]_{_{D^is}}=[a_{(n)}D^is]$,
 then $D^j[a_{(n)}D^is]=[a_{(n)}D^{i+j}s]+
 \sum\alpha_{_t}[a_{_{t}(n_t)}D^{l_t}s]$,
where each $[a_{_{t}(n_t)}D^{l_t}s]$ is a
 normal $S$-word,
$a_{_{t}(n_t)}\bar{s}D^{l_t}<a_{(n)}\bar{s}D^{i+j}$.
\end{enumerate}
\end{lemma}
\noindent{\bf Proof.}
(i)\ If $a=1$, then
$
D^j[s_{(m)}c]=\sum_{t\geq 0}(-1)^{t}\binom{j}{t}\frac{m!}{(m-t)!}s_{(m-t)}D^{j-t}[c].
$
We assume that
$D^{j-t}[c]=[cD^{j-t}]+\sum\beta_{_{t_k}}[w_{_{t_k}}]$ for any $t\geq 0$
by Lemmas \ref{1} and $\ref{l0}$, where $w_{_{t_k}}\in T$ and $w_{_{t_k}}<cD^{j-t}$.
Then
$$
D^j[s_{(m)}c]=s_{(m)}[cD^{j}]+\sum_{t\geq 1}\alpha_{t}s_{(m-t)}[cD^{j-t}]+
\sum_{t\geq 0}\alpha_{t}\beta_{_{t_k}}s_{(m-t)}[w_{_{t_k}}],
$$
where $\alpha_{_t}=(-1)^{t}\binom{j}{t}\frac{m!}{(m-t)!}$.
It is obvious that
$\bar{s}_{(m)}w_{_{0_k}}<\bar{s}_{(m)}cD^{j}$ and
$\bar{s}_{(m-t)}w_{_{t_k}}<\bar{s}_{(m)}cD^{j},\
\bar{s}_{(m-t)}cD^{j-t}<\bar{s}_{(m)}cD^{j}$
for any $1\leq t\leq min\{m,j\}$.

Now let $|a|\geq 1$ and $a=b_{_{1}(n_1)}a_1$ for $b_1\in B$.
Denote for simplicity that $[v]_{s}=[a_{_{1}(n)}s_{(m)}c]$.
Then $[a_{(n)}s_{(m)}c]=b_{_{1}(n_1)}[v]_{s}$ and
\begin{eqnarray*}
D^j(b_{_{1}(n_1)}[v]_{s})&=&\sum_{t\geq 0}(-1)^{t}
\binom{j}{t}\frac{n_1!}{(n_1-t)!}b_{_{1}(n_1-t)}D^{j-t}[v]_{s}\\
&=&b_{_{1}(n_1)}D^{j}[v]_{s}+\sum_{t\geq 1}(-1)^{t}
\binom{j}{t}\frac{n_1!}{(n_1-t)!}b_{_{1}(n_1-t)}D^{j-t}[v]_{s}.
\end{eqnarray*}
By induction on $|a|$,  (i) follows.

(ii)\ If $a=1$, then $D^j(D^is)=D^{i+j}s$
is already a normal $S$-word. Let $|a|\geq 1$.
We may assume that $a=b_{_{1}(n_1)}a_1$ for some
$b_1\in B$. Denote for simplicity that $[u_{_1}]_{_{D^is}}=[a_{_{1}(n)}D^is]$.
Then $[a_{(n)}D^is]=b_{_{1}(n_1)}[u_{_1}]_{_{D^is}}$  and
$$
D^j(b_{_{1}(n_1)}[u_{_1}]_{_{D^is}})=\sum_{t\geq 0}(-1)^{t}
\binom{j}{t}\frac{n_1!}{(n_1-t)!}b_{_{1}(n_1-t)}D^{j-t}[u_{_1}]_{_{D^is}}.
$$
We can get (ii) by induction on $|a|$.
\ \ \ \ $\square$

\begin{definition}
Let $f$ and $g$ be monic polynomials of $C(B,N)$. We have the following compositions.

$\bullet$ If $w=\bar{f}= a_{(n)}\bar{g}_{(m)}c$, $a,c\in T$ and $a, \bar g$ are $D$-free ($a$ maybe empty),
then define
$(f,g)_w=f-[a_{(n)}g_{(m)}c]$,
which is a composition of inclusion.

$\bullet$ If $w=\bar{f}=a_{(n)}\bar{g}D^i$, $a\in T,\ a$ is $D$-free ($a$ maybe empty) and $i\geq 0$,
then define
$(f,g)_w=f-[a_{(n)} D^ig]$,
which is a composition of right inclusion.

$\bullet$ If $w=\bar{f}_{(m)}c=a_{(n)}\bar{g}$,
$a,c\in T$, $a, \bar f$ are $D$-free and $|\bar{f}|+|\bar{g}|>|w|$, then define
$(f,g)_w=[f_{(m)}c]-[a_{(n)}g]$,
which is a composition of intersection.

$\bullet$ If $w=\bar{f}D^i=a_{(n)}\bar{g}$, $a\in T$ is $D$-free and $i>0$,
then define
$(f,g)_w=D^if-[a_{(n)}g]$,
which is a composition of right intersection.

$\bullet$ If $b\in B$ and $n\geq N$,
then $b_{(n)}f$ is referred to as a composition of left multiplication.

$\bullet$ If $b\in B$ and either  $\bar f$ is not D-free and
$0\leq n<N$ or $f$ is not D-free and $n\geq N$,
then $f_{(n)}b$ is referred to as a composition of right multiplication.
\end{definition}

\begin{definition}\label{def1}
Let $S\subset C(B,N)$ be a set of monic polynomials and $h\in C(B,N)$.
Then $h$ is said to be trivial modulo $S$, denoted by
$$
h\equiv0\ mod(S), \ \mbox{ if } \ h=\sum_{i\in I}\alpha_{_i}[u_{_{i}}]_{_{D^{^{l_i}}s_{_i}}},
$$
where each $[u_{_{i}}]_{_{D^{^{l_i}}s_{_i}}}$ is a normal $S$-word, $s_{_i}\in S$
and $\overline{[u_{_{i}}]_{_{D^{^{l_i}}s_{_i}}}}\leq \bar{h}$.

We call $S$  a Gr\"{o}bner-Shirshov basis in $C(B,N)$
if all compositions of polynomials in $S$ are trivial modulo $S$.

If all left (right, resp.) multiplication compositions of elements of $S$ are trivial modulo $S$,
then $S$ is said to be closed under the composition of left (right, resp.) multiplication.
\end{definition}

\begin{definition}\label{def2}
Let $S$ be a monic subset of $C(B,N)$. A normal word $[u]$ is said to
be $S$-irreducible if $u\neq \overline{[v]_{_{D^is}}}$, where $[v]_{_{D^is}}$
is any normal $S$-word.
Denote the set of all $S$-irreducible words by $Irr(S)$.

A  Gr\"{o}bner-Shirshov basis $S$ in $C(B,N)$ is minimal, if there are no
compositions of inclusion and right inclusion of polynomials from $S$, i.e.,
$\bar f\neq \overline{[u]_{D^ig}}$ for any different $f,g\in S$, where
$[u]_{D^ig}$ is any normal $S$-word with $g$.

A Gr\"{o}bner-Shirshov basis in $C(B,N)$
is reduced
provided that
$supp(s)\subseteq  Irr(S\setminus\{s\})$
for every
$s\in S$,
where
$supp(s)=\{u_1,u_2,\dots,u_n\}$
whenever
$s=\sum_{i=1}^n\alpha_i[u_i]$
with
$0\neq\alpha_i\in \mathbf{k}$
and
$u_i\in T$.
In other words,
each
$u_i$
is an
$S\setminus\{s\}$-irreducible word.
\end{definition}

Clearly, a reduced Gr\"{o}bner-Shirshov basis is minimal.

We will prove that for each ideal $I$ of $C(B,N)$,
$I$ has a unique reduced Gr\"{o}bner-Shirshov basis,
see Theorem \ref{t*}.

\begin{remark}\label{remark1}
\begin{enumerate}
\item[(i)] The definition of a Gr\"{o}bner-Shirshov basis in $C(B,N)$
in \cite{BFK04} (see, \cite{BFK04}, Definition 3.11.) is different from in Definition \ref{def1}.

    Comparing to Definition 3.11. in \cite{BFK04}, we have different definitions of the following:
\begin{enumerate}
   \item[1)] normal $S$-word of the first kind.
   \item[2)] composition of inclusion.
   \item[3)] composition of intersection.
    \item[4)] composition of right multiplication.
   \item[5)]composition of polynomials in $S$ is trivial.
\end{enumerate}

\item[(ii)]  Let $S\subset C(B,N)$ be a set of monic polynomials. If $S$ is a
Gr\"{o}bner-Shirshov basis in the sense of \cite{BFK04}, then $S$ is a
Gr\"{o}bner-Shirshov basis in the sense of Definition \ref{def1}, but not conversely.
\end{enumerate}
\end{remark}

The following example comes from \cite{BFK04}, 4.7.4.
\begin{example}\label{ex00}
Let $C(a, N=2\mid a_{(1)}a-a_{(0)}Da=0)$.
Then $S=\{ a_{(1)}a-a_{(0)}Da=0,\ [a_{(0)}a_{(0)}a]=0 \}$  is a reduced
Gr\"{o}bner-Shirshov basis in the sense of Definitions \ref{def1} and \ref{def2}, while $S$  is not a
Gr\"{o}bner-Shirshov basis in the sense of \cite{BFK04}.

In the sense of \cite{BFK04},
 $C(a, N=2\mid a_{(1)}a-a_{(0)}Da=0)$
has a Gr\"{o}bner-Shirshov basis which is not reduced:
\begin{eqnarray*}
S_1&=&\{a_{(1)}a-a_{(0)}Da=0,\ \ [a_{(0)}a_{(0)}a]=0,\ \ [a_{(0)}a_{(1)}a]=0,\\
&&\ \ [a_{(1)}a_{(0)}a]=0,\ \ [a_{(1)}a_{(1)}a]=0\}.
\end{eqnarray*}
\end{example}

\begin{lemma}\label{l7.1}
Let
$[a_{(n)}s_{(m)}c]$ be a normal $S$-word of the first kind.
Then for any $u\in T$ and $p\geq 0$, we have
$$
[a_{(n)}s_{(m)}c]_{(p)}[u]=\sum_t\alpha_{_t}[a_{_{t}(n_t)}s_{(m_t)} c_t],
$$
where each $[a_{_{t}(n_t)}s_{(m_t)} c_{_t}]$ is a normal $S$-word of the first kind
and $a_{_{t}(n_t)}\bar{s}_{(m_t)}c_{_t}\leq
a_{(n)}\bar{s}_{(m)}c^{^{\backslash D}}\natural u$.
In particular, if $p<N$ and $c$ is $D$-free, then
$$
[a_{(n)}s_{(m)}c]_{(p)}[u]=[a_{(n)}s_{(m)}c_{(p)}u]+\sum_t\alpha_{_t}[a_{_{t}(n_t)}s_{(m_t)} c_t],
$$
where  $[a_{_{t}(n_t)}s_{(m_t)} c_{_t}]$ and $[a_{(n)}s_{(m)}c_{(p)}u]$ are normal $S$-words,
$a_{_{t}(n_t)}\bar{s}_{(m_t)}c_{_t}< a_{(n)}\bar{s}_{(m)}c_{(p)}u$.
\end{lemma}
\noindent{\bf Proof.}
Consider first that $a=1$. Then
$$
h:=(s_{(m)}[c])_{(p)}[u]=s_{(m)}([c]_{(p)}[u])+\sum_{t\geq 1}(-1)^t\binom{m}{t}s_{(m-t)}([c]_{(p+t)}[u]).
$$
By Lemma \ref{l0}, $h$ is a linear combination of normal $S$-words
of the first kind.
We may assume that
$
[c]_{(p+t)}[u]=\sum\beta_{_{t_{k}}}[v_{_{t_{k}}}]
$
for any $t\geq 1$, where $v_{_{t_{k}}}\in T$ and
$v_{_{t_{k}}}\leq\overline{[c]_{(p+t)}[u]}\leq c^{^{\backslash D}}\natural u$.
Then the result holds for $a=1$.
In particular, if $p<N$ and $c$ is $D$-free,
then $\overline{[c]_{(p)}[u]}=c_{(p)}u$ by Lemma \ref{1}.
So $ \bar{s}_{(m-t)}v_{_{t_{k}}}<\bar{s}_{(m)}c_{(p)}u$ for $t\geq 1$
and $\bar{h}=\bar{s}_{(m)}c_{(p)}u$.

Now let $|a|\geq1$ and $a=b_{_{1}(n_1)}a_1$ for some
$b_1\in B$. Then
$$
(b_{_{1}(n_1)}[a_{_{1}(n)}s_{(m)}c])_{(p)}[u]
=\sum_{t\geq 0}(-1)^t\binom{n_1}{t}b_{_{1}(n_1-t)}([a_{_{1}(n)}s_{(m)}c]_{(p+t)}[u]).
$$
By induction on $|a|$, we can get the result.
\ \ \ \ $\square$

\begin{lemma}\label{l7.2}
Let
$[a_{(n)}D^is]$ be a normal $S$-word of the second kind.
If $S$ is closed under
the composition of right multiplication, then for
any $u\in T$ and $p\geq 0$, we have
$$
[a_{(n)}D^is]_{(p)}[u]=\sum\alpha_{_t}[v_{_t}]_{_{D^{l_t}s_t}},
$$
where each $[v_{_t}]_{_{D^{l_t}s_t}}$ is a normal $S$-word
and $\overline{[v_{_t}]_{_{D^{l_t}s_t}}}\leq
a_{(n)}\bar{s}^{^{\backslash D}}\natural u$.
In particular, if $i=0,\ p<N$ and $\bar{s}$ is
$D$-free, then
$$
[a_{(n)}s]_{(p)}[u]=[a_{(n)}s_{(p)}u]+\sum\alpha_{_t}[v_{_t}]_{_{D^{l_t}s_t}},
$$
where $\overline{[v_{_t}]_{_{D^{l_t}s_t}}}< a_{(n)}\bar{s}_{(p)}u$,
$[v_{_t}]_{_{D^{l_t}s_t}}$ and $[a_{(n)}s_{(p)}u]$ are normal $S$-words.
\end{lemma}
\noindent{\bf Proof.}
Case 1\ \   $|u|=1,\ a=1,\ i=0$. Suppose that $u=D^jb$
for some $b\in B$. Then
$
s_{(p)}D^jb=\sum_{t\geq 0}\binom{j}{t}\frac{p!}{(p-t)!}D^{j-t}(s_{(p-t)}b).
$
If $s_{(p-t)}b$ is a normal $S$-word, then by Lemma \ref{l4}, $D^{j-t}(s_{(p-t)}b)$ is
a  linear combination of normal $S$-words
$[v_{_k}]_{_{D^{^{i_k}}s_k}}$ and
$\overline{[v_{_k}]_{_{D^{^{i_k}}s_k}}}\leq\overline{D^{j-t}(s_{(p-t)}b)}
=\bar{s}_{(p-t)}D^{j-t}b\leq\bar{s} \natural D^{j}b$.
Let $s_{(p-t)}b$ be not a normal $S$-word and $s_{(p-t)}b\neq 0$.
Since each right multiplication composition of $S$ is trivial,
we suppose that $s_{(p-t)}b=\sum_{k_t}\beta_{k_t}[v_{_{k_t}}]_{D^{^{l_{k_t}}}s_{k_t}}$,
where each $[v_{_{k_t}}]_{D^{^{l_{k_t}}}s_{k_t}}$ is a normal $S$-word and
$\overline{[v_{_{k_t}}]_{D^{^{l_{k_t}}}s_{k_t}}}\leq\overline{s_{(p-t)}b}
\leq
\bar{s}^{^{\backslash D}}\natural b$.
Thus
$
\overline{D^{j-t}([v_{_{k_t}}]_{D^{^{l_{k_t}}}s_{k_t}})}\leq
\overline{D^{j-t}[\bar{s}^{^{\backslash D}}\natural b]}
=\bar{s}^{^{\backslash D}}\natural D^{j-t}b\leq\bar{s}^{^{\backslash D}}\natural D^{j}b.
$
By Lemma \ref{l4}, we see that $s_{(p)}D^jb$ has a presentation
$s_{(p)}D^jb=\sum\beta_{_k}[v_{_k}]_{D^{l_k}s_k}$,
where each $[v_{_k}]_{D^{l_k}s_k}$ is a normal $S$-word,
and $\overline{[v_{_k}]_{D^{l_k}s_k}}\leq\bar{s}^{^{\backslash D}}\natural D^{j}b$.

Case 2\ \   $|u|=1,\ |a|\geq 1,\ i\geq 0$.
We assume that $a=b_{_{1}(n_{1})}\cdots b_{_{k}(n_{k})}b_{_{k+1}}\in T$.
Then we have
\begin{eqnarray*}
h:&=&[a_{(n)}D^is]_{(p)}D^jb\\
&=&\sum\limits_{t_1,\cdots,t_k,t\geq 0}\beta_{t_1\cdots t_kt}
b_{_{1}(n_{1}-t_1)}(\cdots (b_{_{k}(n_{k}-t_k)}(b_{_{k+1}(n-t)}(D^is_{(p+q)}D^jb)))\cdots)\\
&=&\sum\limits_{t_1,\cdots,t_k,t\geq 0}\beta'_{t_1\cdots t_kt}
b_{_{1}(n_{1}-t_1)}(\cdots (b_{_{k}(n_{k}-t_k)}(b_{_{k+1}(n-t)}(s_{(p+q-i)}D^jb)))\cdots),
\end{eqnarray*}
where $q=t_1+\cdots +t_k+t,\
\beta_{t_1\cdots t_kt}=(-1)^{q}\binom{n_1}{t_1}\cdots\binom{n_k}{t_k}\binom{n}{t}$
and $\beta'_{t_1\cdots t_kt}=(-1)^{i}\frac{(p+q)!}{(p+q-i)!}\beta_{t_1\cdots t_kt}$.
Since $n-t\leq N,\ n_l-t_l\leq N$ and $1\leq l\leq k$, we have
$h=\sum\alpha_{_m}[v_{_m}]_{_{D^{l_m}s_m}}$ by Case 1,
where each $[v_{_m}]_{_{D^{l_m}s_m}}$ is a normal $S$-word
and
$$
\overline{[v_{_m}]_{_{D^{l_m}s_m}}}\leq b_{_{1}(n_{1})}(\cdots (b_{_{k}(n_{k})}(b_{_{k+1}(n)}
(\bar{s}^{^{\backslash D}}\natural D^{j}b)))\cdots)=a_{(n)}\bar{s}^{^{\backslash D}}\natural D^{j}b.
$$
In particular, if $i=0,\ p<N$ and $\bar{s}$ is
$D$-free, then $h=[a_{(n)}s_{(p)}D^{j}b]+\sum\alpha_{_m}[v_{_m}]_{_{D^{l_m}s_m}}$,
where $\overline{[v_{_m}]_{_{D^{l_m}s_m}}}<a_{(n)}\bar{s}_{(p)}D^{j}b$.

Case 3\ \  $|u|>1$. We may assume that $u=b_{_{1}(m_1)}u_{_1}\in T$ where $b_1\in B$.
Then
$
[a_{(n)}D^is]_{(p)}[b_{_{1}(m_1)}u_{_1}]=
\sum\limits_{t\geq 0}\binom{p}{t}([a_{(n)}D^is]_{(p-t)}b_{_1})_{(m_1+t)}[u_{_1}].
$
By Case 2, Lemma \ref{l7.1} and induction on $|u|$, we get the result.
\ \ \ \ $\square$

\begin{lemma}\label{l7.3}
Let  $[v]_{_{D^{i}s}}$ be a normal $S$-word
and $S$ be closed under the compositions of left and
right multiplication, Then for any $u\in T$ and $p\geq 0$,
the following statements hold.
\begin{enumerate}
   \item[(i)]\  $[u]_{(p)}[v]_{_{D^{i}s}}=\sum\alpha_{_t}[v_{_t}]_{_{D^{l_t}s_t}}$,
where each $[v_{_t}]_{_{D^{l_t}s_t}}$ is a normal $S$-word
and $\overline{[v_{_t}]_{_{D^{l_t}s_t}}}\leq
u^{^{\backslash D}}\natural\overline{[v]_{_{D^{i}s}}}$.
   \item[(ii)]\  If $p<N$ and $u$ is $D$-free, then
we have
  \begin{enumerate}
  \item[1)]\ If $[v]_{_{D^{i}s}}=[a_{(n)}s_{(m)}c]$, then
\begin{equation}\label{e1.*}
[u]_{(p)}[a_{(n)}s_{(m)}c]=[u_{(p)}a_{(n)}s_{(m)}c]+\sum\alpha_{_t}[v_{_t}]_{_{s_t}},
\end{equation}
where $[v_{_t}]_{_{s_t}}$
and $[u_{(p)}a_{(n)}D^is]$ are normal $S$-words,
$\overline{[v_{_t}]_{_{s_t}}}<u_{(p)}a_{(n)}\bar{s}_{(m)}c$.
  \item[2)]\ If $[v]_{_{D^{i}s}}=[a_{(n)}D^is]$, then
\begin{equation}\label{e1.**}
[u]_{(p)}[a_{(n)}D^is]=[u_{(p)}a_{(n)}D^is]+\sum\alpha_{_t}[u_{_t}]_{_{D^{l_t}s_t}},
\end{equation}
where $[u_{_t}]_{_{D^{l_t}s_t}}$ and $[u_{(p)}a_{(n)}D^is]$ are normal $S$-words,
$\overline{[u_{_t}]_{_{D^{l_t}s_t}}}< u_{(p)}a_{(n)}\bar{s }D^i$.
  \end{enumerate}
\end{enumerate}
\end{lemma}
\noindent{\bf Proof.} (i)\
Case 1\ \  $|u|=1$. Suppose that $u=D^jb$ for some $b\in B$ and $j\geq 0$.
If $j>0$, then
$
D^jb_{(p)}[v]_{_{D^is}}=(-1)^j\frac{(p)!}{(p-j)!}b_{(p-j)}[v]_{_{D^is}}.
$
So it suffices to consider $j=0$.
If $p<N$,
then $b_{(p)}[v]_{_{D^is}}$ is normal. The result is clear.
Now let $p\geq N$.

Case 1.1\ \  Let $[v]_{_{D^{i}s}}=[s_{(m)}c]$ be
a normal $S$-word of the first kind. Then
$
b_{(p)}[s_{(m)}c]=\sum\limits_{t\geq 0}\binom{p}{t}(b_{(p-t)}s)_{(m+t)}[c].
$
If $p-t< N$, then $b_{(p-t)}s$ is a normal $S$-word
of the second kind. By Lemma \ref{l7.2}, we know that
$(b_{(p-t)}s)_{(m+t)}[c]$ can be presented as
a linear combination of normal $S$-words
$[v_{_k}]_{_{D^{^{i_k}}s_k}}$ and
$\overline{[v_{_k}]_{_{D^{^{i_k}}s_k}}}\leq
b_{(p-t)}\bar{s}\natural c\leq b\natural \overline{[s_{(m)}c]}$.
Let $p-t\geq N$. Since $S$ is closed under the composition
of left multiplication, $b_{(p-t)}s$ can be presented
as a linear combination of normal $S$-words $[v_{_l}]_{_{D^{^{i_l}}s_l}}$
and
$\overline{[v_{_l}]_{_{D^{^{i_l}}s_l}}}\leq
\overline{b_{(p-t)}s}\leq b\natural \bar{s}=:w$.
Since $\bar{s}$ is $D$-free, $w$ is $D$-free.
By Lemmas \ref{l7.1} and  \ref{l7.2}, we  get that
$
([v_{_l}]_{_{D^{^{i_l}}s_l}})_{(m+t)}[c]=\sum\alpha_{_k}[u_{_k}]_{_{D^{j_k}s_k}}
$
and $\overline{[u_{_k}]_{_{D^{j_k}s_k}}}\leq
\overline{[v_{_l}]_{_{D^{^{i_l}}s_l}}}^{\backslash D}\natural c
\leq w\natural c=b \natural\overline{[s_{(m)}c]}$.
Hence, the result
holds for $|u|=1$ and $[v]_{_{D^{i}s}}=[s_{(m)}c]$.

Let $[v]_{_{D^is}}=D^is$ be
a normal $S$-word of the second kind. Then
$$
b_{(p)}D^is=\sum_{t\geq 0}\binom{i}{t}\frac{(p)!}{(p-t)!}D^{i-t}(b_{(p-t)}s).
$$
If $p-t< N$, then $b_{(p-t)}s$ is a normal $S$-word.
Let $p-t\geq N$. Since $S$ is closed under the composition
of left multiplication, $b_{(p-t)}s$ can be presented
as a linear combination of normal $S$-words $[v_{_l}]_{_{D^{^{i_l}}s_l}}$
and
$\overline{[v_{_l}]_{_{D^{^{i_l}}s_l}}}\leq
\overline{b_{(p-t)}s}\leq b\natural \bar{s}$.
By Lemma \ref{1}, we have
$\overline{D^{i-t}([v_{_l}]_{_{D^{^{i_l}}s_l}})}\leq b\natural \bar{s}D^{i}$.
By Lemma \ref{l4}, the result holds for $|u|=1$ and
$[v]_{_{D^is}}=D^is$.

Case 1.2\ \  Let $a=b_{_{1}(n_1)}a_{_1}$
and $[v]_{_{D^is}}=[a_{(n)}s_{(m)}c]$ or
$[v]_{_{D^is}}=[a_{(n)}D^is]$. Then we can denote
$[v]_{_{D^is}}=b_{_{1}(n_1)}[v']_{_{D^is}}$,
where $[v']_{_{D^is}}=[a_{_{1}(n)}s_{(m)}c]$ or
$[v']_{_{D^is}}=[a_{_{1}(n)}D^is]$.
Then
\begin{equation}\label{e1.***}
b_{(p)}[v]_{_{D^is}}=b_{(p)}(b_{_{1}(n_1)}[v']_{_{D^is}})
=\sum_{t\geq 1}(-1)^t\binom{p}{t}b_{(p-t)}(b_{1(n_1+t)}[v']_{_{D^is}}).
\end{equation}
Since $|a_1|<|a|$, by induction on $|a|$, we have
$
b_{1(n_1+t)}[v']_{_{D^is}}=\sum\beta_{_k}[w_{_k}]_{_{D^{j_k}s_k}}
$
and $\overline{[w_{_k}]_{_{D^{j_k}s_k}}}\leq b_{1}\natural \overline{[v']_{_{D^is}}}$.
If $p-t\geq N$, we can make $p-t$ smaller by using the formula (\ref{e1.***}) again.
So we may assume $0\leq p-t<N$. Thus each $b_{(p-t)}[w_{_k}]_{_{D^{j_k}s_k}}$
is a normal $S$-word and
$\overline{b_{(p-t)}[w_{_k}]_{_{D^{j_k}s_k}}}
\leq b_{(N-1)}b_{1}\natural
\overline{[v']_{_{D^is}}}
=b\natural\overline{[v]_{_{D^is}}}$.
Therefore,  the result holds for $|u|=1$.

Case 2\ \  $|u|>1$. We may assume that
$u=b_{_{1}(n_1)}u_{_1}\in T$ for some $b_1\in B$.
Then
\begin{eqnarray*}
[b_{_{1}(n_1)}u_{_1}]_{(p)}[v]_{_{D^is}}
&=&\sum_{t\geq 0}
(-1)^t\binom{n_1}{t}b_{_{1}(n_1-t)}([u_{_1}]_{(p+t)}[v]_{_{D^is}})\\
&=&b_{_{1}(n_1)}([u_{_1}]_{(p)}[v]_{_{D^is}})+\sum_{t\geq 1}
(-1)^t\binom{n_1}{t}b_{_{1}(n_1-t)}([u_{_1}]_{(p+t)}[v]_{_{D^is}}).
\end{eqnarray*}
Since $n_1-t<N$, we have the result by induction on $|u|$.

(ii)\ If $p<N$ and $u$ is $D$-free, by associative identity and
induction on $|u|$, we can get identities (\ref{e1.*}) and (\ref{e1.**}).
\ \ \ \ $\square$

\begin{lemma}\label{l8}
Let $S$  be closed under
the compositions of left and right multiplication.
Then any $S$-word can be presented as a linear combination
of normal $S$-words.
\end{lemma}
\noindent{\bf Proof.}
Let $(u)_{_{D^{i}s}}$ be an $S$-word.
Induction on $|u|_{_{D^{i}s}}$.
The result holds trivially if $|u|_{_{D^{i}s}}=1$.
Let $|u|_{_{D^{i}s}}>1$.
There are essentially two possibilities:
$$
(u)_{_{D^is}}=(a)_{(n)}(u_{_1})_{_{D^is}} \ \mbox{or} \ (u)_{_{D^is}}=(u_{_1})_{_{D^is}{(n)}}(a),
$$
where $(a)$ is a word, $(u_{_1})_{_{D^is}}$ is an $S$-word
and $|u_{_1}|_{_{D^{i}s}}<|u|_{_{D^{i}s}}$. By Lemma \ref{l0}, we may assume that $(a)=[a]$ is a normal word.
By induction, $(u_{_1})_{_{D^is}}$ is a linear combination
of normal $S$-words.
So it suffices to consider
$[a]_{(n)}[u_{_1}]_{_{D^is}}$ and $[u_{_1}]_{_{D^is}{(n)}}[a]$.
By Lemmas \ref{l7.1}, \ref{l7.2} and \ref{l7.3},
we complete the proof.
\ \  \ \ $\square$

\ \

The following lemma plays a key role in proving Theorem \ref{t1}.

\begin{lemma}\label{ll}
Let $S$ be a Gr\"{o}bner-Shirshov basis in $C(B,N)$ and $s_{_1},s_{_2}\in S$.
Let $[u_{_1}]_{D^{^{j_{_1}}}s_{_1}}$ and $[u_{_2}]_{D^{^{j_{_2}}}s_{_2}}$ be normal $S$-words.
If $w=\overline{[u_{_1}]_{D^{^{j_{_1}}}s_{_1}}}=\overline{[u_{_2}]_{D^{^{j_{_2}}}s_{_2}}}$,
then
\begin{equation}\label{*}
h:=[u_{_2}]_{D^{^{j_{_2}}}s_{_2}}-[u_{_1}]_{D^{^{j_{_1}}}s_{_1}}
=\sum_t\beta_{_t}[v_{_{t}}]_{_{D^{l_t}s_{_{t}}}},
\end{equation}
where each $[v_{_{t}}]_{_{D^{l_t}s_{_{t}}}}$ is a normal $S$-word and
$\overline{[v_{_{t}}]_{_{D^{l_t}s_{_{t}}}}}<w$.
\end{lemma}
\noindent{\bf Proof.}
There are three cases to consider:

Case 1\ \  $[u_{_1}]_{D^{^{j_{_1}}}s_{_1}}=[a_{_{1}(n_1)}s_{_{1}(m_1)}c_{_1}]$
and $[u_{_2}]_{D^{^{j_{_2}}}s_{_2}}=[a_{_{2}(n_2)}s_{_{2}(m_2)}c_{_2}]$.
Then
$w=a_{_{1}(n_1)}\overline{s_{_1}}_{(m_1)}c_{_1}
=a_{_{2}(n_2)}\overline{s_{_2}}_{(m_2)}c_{_2}$. We have
three subcases to discuss.

Case 1.1\ \  $\overline{s_{_1}}$ and $\overline{s_{_2}}$
are mutually disjoint.
We assume that $a_{_2}=a_{_{1}(n_1)}\overline{s_{_1}}_{(m_1)}a$
and $c_{_1}=a_{(n_2)}\overline{s_{_2}}_{(m_2)}c_{_2}$,
where $a\in T$ or $a=1$. Suppose that
$s_i=[\overline{s_i}]+\sum\alpha_{_{k_i}}[u_{_{k_i}}]$,
where $u_{_{k_i}}\in T,\ u_{_{k_i}}<\overline{s_i},\ i\in \{1,2\}$.
So $u_{_{k_i}}^{\backslash D}<\overline{s_i}$.
By Lemmas \ref{l0} and \ref{1}, we have
$[\overline{s_{_2}}]_{(m_2)}[c_{_2}]=
[\overline{s_{_2}}_{(m_2)}c_{_2}]+\sum\gamma_{_t}[v_{_t}]$,
where $v_{_t}\in T$ and $v_{_t}<\overline{s_{_2}}_{(m_2)}c_{_2}$.
Then
\begin{eqnarray*}
&&[a_{_{1}(n_1)}s_{_{1}(m_1)}c_{_1}]\\
&=&[a_{_{1}(n_1)}s_{_{1}(m_1)}a_{(n_2)}\overline{s_{_2}}_{(m_2)}c_{_2}]\\
&=&[a_{_{1}(n_1)}s_{_{1}(m_1)}a_{(n_2)}[\overline{s_{_2}}]_{(m_2)}[c_{_2}]]
+\sum\gamma_{_t}[a_{_{1}(n_1)}s_{_{1}(m_1)}a_{(n_2)}v_{_t}]\\
&=&[a_{_{1}(n_1)}s_{_{1}(m_1)}a_{(n_2)}s_{_{2}(m_2)}[c_{_2}]]-
[a_{_{1}(n_1)}s_{_{1}(m_1)}a_{(n_2)}(s_{2}-[\overline{s_{_2}}])_{(m_2)}[c_{_2}]]\\
&&+\sum\gamma_{_t}[a_{_{1}(n_1)}s_{_{1}(m_1)}a_{(n_2)}v_{_t}]\\
&=&[a_{_{1}(n_1)}s_{_{1}(m_1)}a_{(n_2)}s_{_{2}(m_2)}[c_{_2}]]-
\sum\alpha_{_{k_2}}[a_{_{1}(n_1)}s_{_{1}(m_1)}a_{(n_2)}[u_{_{k_2}}]_{(m_2)}[c_{_2}]]\\
&&+\sum\gamma_{_t}[a_{_{1}(n_1)}s_{_{1}(m_1)}a_{(n_2)}v_{_t}]\\
&=:&A.
\end{eqnarray*}
By Lemmas \ref{l2} and \ref{l0},
we have
\begin{eqnarray*}
\overline{[a_{_{1}(n_1)}s_{_{1}(m_1)}a_{(n_2)}[u_{_{k_2}}]_{(m_2)}[c_{_2}]]}
&=&a_{_{1}(n_1)}\bar{s}_{_{1}(m_1)}a_{(n_2)}\overline{[u_{_{k_2}}]_{(m_2)}[c_{_2}]}\\
&\leq & a_{_{1}(n_1)}\bar{s}_{_{1}(m_1)}a_{(n_2)}u_{_{k_2}}^{\backslash D}\natural c_{_2}\\
&<&a_{_{1}(n_1)}\bar{s}_{_{1}(m_1)}a_{(n_2)}\overline{s_{_2}}_{(m_2)}c_{_2}=w,\\
\overline{[a_{_{1}(n_1)}s_{_{1}(m_1)}a_{(n_2)}v_{_t}]}
&=&a_{_{1}(n_1)}\bar{s}_{_{1}(m_1)}a_{(n_2)}v_{_t}\\
&<&a_{_{1}(n_1)}\bar{s}_{_{1}(m_1)}a_{(n_2)}\overline{s_{_2}}_{(m_2)}c_{_2}=w.
\end{eqnarray*}
By Lemma \ref{l7.3}, we have
\begin{eqnarray*}
&&[a_{_{2}(n_2)}s_{_{2}(m_2)}c_{_2}]\\
&=&[a_{_{1}(n_1)}\overline{s_{_1}}_{(m_1)}a_{(n_2)}s_{_{2}(m_2)}c_{_2}]\\
&=&[a_{_{1}(n_1)}[\overline{s_{_1}}]_{(m_1)}[a_{(n_2)}s_{_{2}(m_2)}c_{_2}]]+
\sum\beta_{_l}[q_{_{l}}]_{_{D^{k_l}s_{_{l}}}}\\
&=&[a_{_{1}(n_1)}s_{_{1}(m_1)}a_{(n_2)}s_{_{2}(m_2)}c_{_2}]-
[a_{_{1}(n_1)}(s_1-[\overline{s_{_1}}])_{(m_1)}[a_{(n_2)}s_{_{2}(m_2)}c_{_2}]]\\
&&+\sum\beta_{_l}[q_{_{l}}]_{_{D^{k_l}s_{_{l}}}}\\
&=&[a_{_{1}(n_1)}s_{_{1}(m_1)}a_{(n_2)}s_{_{2}(m_2)}c_{_2}]-
\sum\alpha_{_{k_1}}[a_{_{1}(n_1)}[u_{_{k_1}}]_{(m_1)}[a_{(n_2)}s_{_{2}(m_2)}c_{_2}]]\\
&&+\sum\beta_{_l}[q_{_{l}}]_{_{D^{k_l}s_{_{l}}}}\\
&=:&C,
\end{eqnarray*}
where $\overline{[q_{_{l}}]_{_{D^{k_l}s_{_{l}}}}}<w$.
Then by Lemma \ref{l7.3} and Corollary \ref{c1}, we have
\begin{eqnarray*}
\overline{[a_{_{1}(n_1)}[u_{_{k_1}}]_{(m_1)}[a_{(n_2)}s_{_{2}(m_2)}c_{_2}]]}
&=&a_{_{1}(n_1)}\overline{[u_{_{k_1}}]_{(m_1)}[a_{(n_2)}s_{_{2}(m_2)}c_{_2}]}\\
&\leq &a_{_{1}(n_1)}u_{_{k_1}}^{\backslash D}\natural \overline{[a_{(n_2)}s_{_{2}(m_2)}c_{_2}]}\\
&<&a_{_{1}(n_1)}\overline{s_{_1}}_{(m_1)}a_{(n_2)}\overline{s_{_2}}_{(m_2)}c_{_2}=w.
\end{eqnarray*}
Therefore,
\begin{eqnarray*}
h&=&C-A\\
&=&\sum\alpha_{_{k_2}}[a_{_{1}(n_1)}s_{_{1}(m_1)}a_{(n_2)}[u_{_{k_2}}]_{(m_2)}[c_{_2}]]-
\sum\gamma_{_t}[a_{_{1}(n_1)}s_{_{1}(m_1)}a_{(n_2)}v_{_t}]\\
&&-\sum\alpha_{_{k_1}}[a_{_{1}(n_1)}[u_{_{k_1}}]_{(m_1)}[a_{(n_2)}s_{_{2}(m_2)}c_{_2}]]
+\sum\beta_{_l}[q_{_{l}}]_{_{D^{k_l}s_{_{l}}}}.
\end{eqnarray*}
This shows that the result holds for Case 1.1.

Case 1.2\ \  One of $\overline{s_{_1}}$ and $\overline{s_{_2}}$
is a subword of the other.
We assume that $\overline{s_{_2}}$ is a subword of $\overline{s_{_1}}$, i.e.,
$w'_1=\overline{s_{_1}}=a_{(n_2)}\overline{s_{_2}}_{(m_2)}c$,
$a_{_2}=a_{_{1}(n_1)}a$ and $c_{_2}=c_{(m_1)}c_{_1}$,
where $a,c\in T$ are D-free or one of them is empty.
Then by Lemma \ref{l7.2},
\begin{eqnarray*}
h&=&[a_{_{2}(n_2)}s_{_{2}(m_2)}c_{_2}]-[a_{_{1}(n_1)}s_{_{1}(m_1)}c_{_1}]\\
&=&[a_{_{1}(n_1)}a_{(n_2)}s_{_{2}(m_2)}c_{(m_1)}c_{_1}]-[a_{_{1}(n_1)}s_{_{1}(m_1)}c_{_1}]\\
&=&
[a_{_{1}(n_1)}[a_{(n_2)}s_{_{2}(m_2)}c]_{(m_1)}[c_{_1}]]-[a_{_{1}(n_1)}s_{_{1}(m_1)}c_{_1}]
+\sum\gamma_{_t}[q_{_{t}}]_{_{D^{k_t}s_{_{t}}}}\\
&=&
-[a_{_{1}(n_1)}(s_{_1},s_{_2})_{_{w'_1}(m_1)}[c_{_1}]]
+\sum\gamma_{_t}[q_{_{t}}]_{_{D^{k_t}s_{_{t}}}},
\end{eqnarray*}
where $(s_{_1},s_{_2})_{w'_1}=s_{_1}-[a_{(n_2)}s_{_{2}(m_2)}c]$,
each $[q_{_{t}}]_{_{D^{k_t}s_{_{t}}}}$ is a normal $S$-word and
$\overline{[q_{_{t}}]_{_{D^{k_t}s_{_{t}}}}}<w$.
Since $S$ is a Gr\"{o}bner-Shirshov basis,
we know that $(s_{_1},s_{_2})_{w'_1}$ is a linear
combination of normal $S$-words $[v'_{_t}]_{D^{^{l'_t}}s_t}$
and $\overline{[v'_{_t}]_{D^{^{l'_t}}s_t}}\leq\overline{(s_{_1},s_{_2})_{w'_1}}<w'_1$.
Since $w'_1$ is $D$-free, by Corollary \ref{c1},
$\overline{[v'_{_t}]_{D^{^{l'_t}}s_t}}^{^{\backslash D}}\natural c_{_1}
\leq\overline{(s_{_1},s_{_2})_{w'_1}}^{^{\backslash D}}\natural c_{_1}<w'_{_{1}(m_1)}c_{_1}$.
By Lemmas \ref{l7.1}, \ref{l7.2} and \ref{l7.3},
we get that the result holds for Case 1.2.

Case 1.3\ \  $\overline{s_{_1}}$ and $\overline{s_{_2}}$
have a nonempty intersection as subword of $w$,
$\overline{s_{_1}}$ and $\overline{s_{_2}}$ are
not subwords of each other.
Assume that $a_{_2}=a_{_{1}(n_1)}a,\ c_{_1}=c_{(m_2)}c_{_2}$ and
$\overline{s_{_1}}_{(m_1)}c=a_{(n_2)}\overline{s_{_2}}=w'_1$,
where $a,c\in T$ are D-free and
$|\overline{s_{_1}}|+|\overline{s_{_2}}|>|w'_1|$.
Then by Lemmas \ref{l7.1} and \ref{l7.2},
we have
\begin{eqnarray*}
h&=&[a_{_{2}(n_2)}s_{_{2}(m_2)}c_{_2}]-[a_{_{1}(n_1)}s_{_{1}(m_1)}c_{_1}]\\
&=&[a_{_{1}(n_1)}a_{(n_2)}s_{_{2}(m_2)}c_{_2}]-[a_{_{1}(n_1)}s_{_{1}(m_1)}c_{(m_2)}c_{_2}]\\
&=&
[a_{_{1}(n_1)}[a_{(n_2)}s_{_2}]_{(m_2)}[c_{_2}]]-[a_{_{1}(n_1)}[s_{_{1}(m_1)}c]_{(m_2)}[c_{_2}]]
+\sum\gamma_{_t}[q_{_{t}}]_{_{D^{k_t}s_{_{t}}}}\\
&=&
-[a_{_{1}(n_1)}(s_{_1},s_{_2})_{_{w'_1}(m_2)}[c_{_2}]]
+\sum\gamma_{_t}[q_{_{t}}]_{_{D^{k_t}s_{_{t}}}},
\end{eqnarray*}
where $(s_{_1},s_{_2})_{w'_1}=[s_{_{1}(m_1)}c]-[a_{(n_2)}s_{_2}]$,
each $[q_{_{t}}]_{_{D^{k_t}s_{_{t}}}}$ is a normal $S$-word and
$\overline{[q_{_{t}}]_{_{D^{k_t}s_{_{t}}}}}<w$.
Similar to Case 1.2, we get that the result holds for Case 1.3.

Case 2\ \ $[u_{_1}]_{D^{^{j_{_1}}}s_{_1}}=[a_{_{1}(n_1)}s_{_{1}(m_1)}c_{_1}]$
and $[u_{_2}]_{D^{^{j_{_2}}}s_{_2}}=[a_{_{2}(n_2)}D^{^{j_{_2}}}s_{_{2}}]$.
Then $w=a_{_{1}(n_1)}\overline{s_{_1}}_{(m_1)}c_{_1}
=a_{_{2}(n_2)}\overline{s_{_2}}D^{j_{_2}}$.
We also have three subcases to consider.

Case 2.1\ \   $\overline{s_{_1}}$ and $\overline{s_{_2}}$
are mutually disjoint.
We assume that $a_{_2}=a_{_{1}(n_1)}\overline{s_{_1}}_{(m_1)}a$
and $c_{_1}=a_{(n_2)}\overline{s_{_2}}D^{j_{_2}}$,
where $a$ may be empty. By Lemmas \ref{1} and \ref{l7.3}, we have
\begin{eqnarray*}
h&=&[a_{_{2}(n_2)}D^{j_{_2}}s_{_2}]-[a_{_{1}(n_1)}s_{_{1}(m_1)}c_{_1}]\\
&=&[a_{_{1}(n_1)}\overline{s_{_1}}_{(m_1)}a_{(n_2)}D^{j_{_2}}s_{_2}]
-[a_{_{1}(n_1)}s_{_{1}(m_1)}a_{(n_2)}\overline{s_{_2}}D^{j_{_2}}]\\
&=&
[a_{_{1}(n_1)}[\overline{s_{_1}}]_{(m_1)}a_{(n_2)}D^{j_{_2}}s_{_2}]
-[a_{_{1}(n_1)}s_{_{1}(m_1)}a_{(n_2)}D^{j_{_2}}[\overline{s_{_2}}]]
+\sum\gamma_{_t}[q_{_{t}}]_{_{D^{k_t}s_{_{t}}}}\\
&=&
[a_{_{1}(n_1)}s_{_{1}(m_1)}a_{(n_2)}D^{j_{_2}}(s_{_2}-[\overline{s_{_2}}])]-
[a_{_{1}(n_1)}(s_{_1}-[\overline{s_{_1}}])_{(m_1)}a_{(n_2)}D^{j_{_2}}s_{_2}]\\
&&+\sum\gamma_{_t}[q_{_{t}}]_{_{D^{k_t}s_{_{t}}}},
\end{eqnarray*}
where each $[q_{_{t}}]_{_{D^{k_t}s_{_{t}}}}$ is a normal $S$-word and
$\overline{[q_{_{t}}]_{_{D^{k_t}s_{_{t}}}}}<w$.
Similar to Case 1.1, we get the result.

Case 2.2\ \   $\overline{s_{_1}}$ is a subword of $\overline{s_{_2}}$.
Let $w'_1=\overline{s_{_2}}=a_{(n_1)}\overline{s_{_1}}_{(m_1)}c$,
$a_1=a_{_{2}(n_2)}a$ and $c_{_1}=cD^{j_{_2}}$, where $a,c\in T$.
Then by Lemma \ref{l4}, we have
\begin{eqnarray*}
h&=&[a_{_{2}(n_2)}D^{j_{_2}}s_{_2}]-[a_{_{1}(n_1)}s_{_{1}(m_1)}c_{_1}]\\
&=&[a_{_{2}(n_2)}D^{j_{_2}}s_{_2}]-[a_{_{2}(n_2)}a_{(n_1)}s_{_{1}(m_1)}cD^{j_{_2}}]\\
&=&
[a_{_{2}(n_2)}D^{j_{_2}}s_{_2}]-[a_{_{2}(n_2)}D^{j_{_2}}[a_{(n_1)}s_{_{1}(m_1)}c]]
+\sum\gamma_{_t}[q_{_{t}}]_{_{D^{k_t}s_{_{t}}}}\\
&=&[a_{_{2}(n_2)}D^{j_{_2}}((s_{_2},s_{_1})_{w'_1})]
+\sum\gamma_{_t}[q_{_{t}}]_{_{D^{k_t}s_{_{t}}}},
\end{eqnarray*}
where $(s_{_2},s_{_1})_{w'_1}=s_{_2}-[a_{(n_1)}s_{_{1}(m_1)}c]$,
each $[q_{_{t}}]_{_{D^{k_t}s_{_{t}}}}$ is a normal $S$-word and
$\overline{[q_{_{t}}]_{_{D^{k_t}s_{_{t}}}}}<w$.
Now the result follows from Lemma \ref{l4}.

Case 2.3\ \   $\overline{s_{_1}}$ and $\overline{s_{_2}}$ have
a nonempty intersection,
$\overline{s_{_1}}$ and $\overline{s_{_2}}$ are
not subwords of each other.
Assume that $w'_1=\overline{s_{_1}}_{(m_1)}c=a_{(n_2)}\overline{s_{_2}}$,
$a_{_2}=a_{_{1}(n_1)}a$ and $c_{_1}=cD^{j_{_2}}$, where
$a,c\in T$ and $|\overline{s_{_1}}|+|\overline{s_{_2}}|>|w'_1|$.
Then
\begin{eqnarray*}
h&=&[a_{_{2}(n_2)}D^{j_{_2}}s_{_2}]-[a_{_{1}(n_1)}s_{_{1}(m_1)}c_{_1}]\\
&=&[a_{_{1}(n_1)}a_{(n_2)}D^{j_{_2}}s_{_2}]-[a_{_{1}(n_1)}s_{_{1}(m_1)}cD^{j_{_2}}]\\
&=&
[a_{_{1}(n_1)}D^{j_{_2}}[a_{(n_2)}s_{_2}]]-[a_{_{1}(n_1)}D^{j_{_2}}[s_{_{1}(m_1)}c]]
+\sum\gamma_{_t}[q_{_{t}}]_{_{D^{k_t}s_{_{t}}}}\\
&=&-[a_{_{1}(n_1)}D^{j_{_2}}((s_{_1},s_{_2})_{w'_1})]
+\sum\gamma_{_t}[q_{_{t}}]_{_{D^{k_t}s_{_{t}}}},
\end{eqnarray*}
where $(s_{_1},s_{_2})_{w'_1}=[s_{_{1}(m_1)}a]-[a_{(n_2)}s_{_2}]$,
each $[q_{_{t}}]_{_{D^{k_t}s_{_{t}}}}$ is a normal $S$-word and
$\overline{[q_{_{t}}]_{_{D^{k_t}s_{_{t}}}}}<w$.
By Lemma \ref{l4}, we can get the result.

Case 3\ \   $[u_{_1}]_{D^{^{j_{_1}}}s_{_1}}=[a_{_{1}(n_1)}D^{^{j_{_1}}}s_{_{1}}]$
and $[u_{_2}]_{D^{^{j_{_2}}}s_{_2}}=[a_{_{2}(n_2)}D^{^{j_{_2}}}s_{_{2}}]$.
Then  $w_1=a_{_{1}(n_1)}\overline{s_{_1}}D^{j_{_1}}
=a_{_{2}(n_2)}\overline{s_{_2}}D^{j_{_2}}$.
We may assume that $j_{_2}\geq j_{_1}$ and put $i=j_{_2}-j_{_1}\geq 0$.
Denote $w_1=w'_1D^{j_{_1}}$, where
$w'_1=a_{_{1}(n_1)}\overline{s_{_1}}=a_{_{2}(n_2)}\overline{s_{_2}}D^i$.
Then there are two possibilities.

Case 3.1\ \  $\overline{s_{_2}}$ is a subword of $\overline{s_{_1}}$,
i.e., $a_{_2}=a_{_{1}(n_1)}a$ and $w'=\overline{s_{_1}}=a_{(n_2)}\overline{s_{_2}}D^i$,
where $a\in T$ is D-free.
Then
\begin{eqnarray*}
h&=&[a_{_{2}(n_2)}D^{j_{_2}}s_{_2}]-[a_{_{1}(n_1)}D^{j_{_1}}s_{_1}]\\
&=&[a_{_{1}(n_1)}a_{(n_2)}D^{j_{_1}+i}s_{_2}]-[a_{_{1}(n_1)}D^{j_{_1}}s_{_1}]\\
&=&
[a_{_{1}(n_1)}D^{j_{_1}}[a_{(n_2)}D^{i}s_{_2}]-[a_{_{1}(n_1)}D^{j_{_1}}s_{_1}]
+\sum\gamma_{_t}[q_{_{t}}]_{_{D^{k_t}s_{_{t}}}}\\
&=&-[a_{_{1}(n_1)}D^{j_{_1}}((s_{_1},s_{_2})_{w'})]
+\sum\gamma_{_t}[q_{_{t}}]_{_{D^{k_t}s_{_{t}}}},
\end{eqnarray*}
where $(s_{_1},s_{_2})_{w'}=s_{_1}-[a_{(n_2)}D^is_{_2}]$,
each $[q_{_{t}}]_{_{D^{k_t}s_{_{t}}}}$ is a normal $S$-word and
$\overline{[q_{_{t}}]_{_{D^{k_t}s_{_{t}}}}}<w$.
By Lemma \ref{l4}, we can get the result.

Case 3.2\ \   $\overline{s_{_1}}$ is a subword of $\overline{s_{_2}}D^i$,
i.e., $a_1=a_{_{2}(n_2)}a$ and $w'=\overline{s_{_2}}D^i=a_{(n_1)}\overline{s_{_1}}$,
where $a\in T$ is D-free.
Then
\begin{eqnarray*}
h&=&[a_{_{2}(n_2)}D^{j_{_2}}s_{_2}]-[a_{_{1}(n_1)}D^{j_{_1}}s_{_1}]\\
&=&[a_{_{2}(n_2)}D^{j_{_1}+i}s_{_2}]-[a_{_{2}(n_2)}a_{(n_1)}D^{j_{_1}}s_{_1}]\\
&=&
[a_{_{2}(n_2)}D^{j_{_1}+i}s_{_2}]-[a_{_{2}(n_2)}D^{j_{_1}}[a_{(n_1)}s_{_1}]]
+\sum\gamma_{_t}[q_{_{t}}]_{_{D^{k_t}s_{_{t}}}}\\
&=&
[a_{_{2}(n_2)}D^{j_{_1}}((s_{_2},s_{_1})_{w'})]
+\sum\gamma_{_t}[q_{_{t}}]_{_{D^{k_t}s_{_{t}}}},
\end{eqnarray*}
where $(s_{_2},s_{_1})_{w'}=D^is_{_2}-[a_{(n_1)}s_{_1}]$,
each $[q_{_{t}}]_{_{D^{k_t}s_{_{t}}}}$ is a normal $S$-word and
$\overline{[q_{_{t}}]_{_{D^{k_t}s_{_{t}}}}}<w$.
By Lemma \ref{l4}, we can get the result.

This shows (\ref{*}).
\ \ \ \  $\square$

\begin{lemma}\label{l15}
Let $S$ be a monic subset of $C(B,N)$. Then for any $f\in C(B,N)$,
$$
f=\sum_j\alpha_{_j}[u_{_{j}}]_{_{D^{l_j}s_{_j}}}
+\sum_k\alpha_{_k}[u_{_k}],
$$
where  each $[u_{_{j}}]_{_{D^{l_j}s_{_j}}}$ is a normal $S$-word,
$\overline{[u_{_{j}}]_{_{D^{l_j}s_{_j}}}}\leq\bar f,\ [u_{_k}]\in Irr(S),\
u_{_k}\leq \bar f,\
\alpha_{_j},\alpha_{_k}\in {\bf k}$ and $s_{_j}\in S$.
\end{lemma}
\noindent{\bf Proof.}
Induction on $\bar{f}$.
Let $f=\sum\alpha_{_i}[w_{_i}]$, where $w_{_i}\in T$ and $\alpha_{_i}\neq 0$.
 We may assume without loss of generality
that $w_{_1}>w_{_2}>\cdots$. If $[w_{_1}]\in Irr(S)$,
then let $f_1=f-\alpha_{_1}[w_{_1}]$.
If $w_{_1}=\overline{[v]_{_{D^is}}}$ for some normal $S$-word $[v]_{_{D^is}}$,
then let $f_1=f-\alpha_{_1}[v]_{_{D^is}}$.
In all cases, we have $\bar{f}>\overline{f_1}$. Then the result follows
from the induction on $\bar{f}$.
\ \ \ \  $\square$

\ \

 The following theorem is the main result in this paper.

\begin{theorem} (Composition-Diamond lemma for associative conformal algebras)\label{t1}
Let $C(B,N)$ be the free associative conformal algebra
generated by a well-ordered set $B$ with a bounded locality $N$.
Let $S$ be a monic subset of $C(B,N)$ and $Id(S)$ be
the ideal of $C(B,N)$ generated by $S$.
Then the following statements are equivalent.
\begin{enumerate}
 \item[(i)]\ $S$ is a Gr\"{o}bner-Shirshov basis in $C(B,N)$.
 \item[(ii)]\ If $0\neq f\in Id(S)$, then $f=\sum\alpha\overline{[v]_{_{D^is}}}$
 for some normal $S$-word $[v]_{_{D^is}}$.
  \item[$(ii)'$]\ If $0\neq f\in Id(S)$, then
  $f=\sum\alpha_{_i}[u_{_{i}}]_{_{D^{j_i}s_{_{i}}}}$, where  each
$\alpha_i\in \mathbf{k},\ [u_{_{i}}]_{_{D^{j_i}s_{_{i}}}}$ is a normal $S$-word and $ \bar f=\overline{[u_{_{1}}]_{_{D^{j_1}s_{_{1}}}}}>\overline{[u_{_{2}}]_{_{D^{j_2}s_{_{2}}}}}>\cdots$.
 \item[(iii)]\ $Irr(S)$ is a $\mathbf{k}$-basis of $C(B,N|S)=C(B,N)/Id(S)$.
\end{enumerate}
\end{theorem}
\noindent{\bf Proof.}
$(i)\Rightarrow (ii)$. Let $S$ be a Gr\"{o}bner-Shirshov basis
and $0\neq f\in Id(S)$. By Lemma \ref{l8}, $f$ has an expression
$f=\sum\alpha_{_i}[u_{_{i}}]_{_{D^{j_i}s_{_{i}}}}$,
where each
$[u_{_{i}}]_{_{D^{j_i}s_{_{i}}}}$ is a normal $S$-word.
Denote $w_i=\overline{[u_{_{i}}]_{_{D^{j_i}s_{_{i}}}}}$.
We may assume that
$$
w=w_1=w_2=\cdots =w_l> w_{l+1}\geq w_{l+2}\geq \cdots
$$
\noindent for some $l\geq 1$.

If $l=1$ then the result is clear.
Suppose that $l>1$. Then by Lemma \ref{ll}, we have
\begin{eqnarray*}
\alpha_{_1}[u_{_{1}}]_{_{D^{j_1}s_{_{1}}}}+\alpha_{_2}[u_{_{2}}]_{_{D^{j_2}s_{_{2}}}}
&=&(\alpha_{_1}+\alpha_{_2})[u_{_{1}}]_{_{D^{j_1}s_{_{1}}}}
+\alpha_{_2}([u_{_{2}}]_{_{D^{j_2}s_{_{2}}}}-[u_{_{1}}]_{_{D^{j_1}s_{_{1}}}})\\
&=&(\alpha_{_1}+\alpha_{_2})[u_{_{1}}]_{_{D^{j_1}s_{_{1}}}}
+\sum\alpha_{_2}\beta_{_t}[v_{_{t}}]_{_{D^{l_t}s_{_{t}}}},
\end{eqnarray*}
where each $[v_{_{t}}]_{_{D^{l_t}s_{_{t}}}}$ is a normal $S$-word and
$\overline{[v_{_{t}}]_{_{D^{l_t}s_{_{t}}}}}<w$.

Then (ii) follows from induction on $(w,\ l)$.

$(ii)\Rightarrow (ii)'$. This part is clear.

$(ii)'\Rightarrow (iii)$.
 The set $Irr(S)$ generates the algebra $C(B,N)/Id(S)$ as vector space
 by Lemma \ref{l15}.
On the other hand, suppose that $h=\sum\alpha_{_t}[u_{_t}]=0$ in $C(B,N)/Id(S)$,
where each $[u_{_t}]\in{Irr(S)}$. This means that $h\in{Id(S)}$.
Then each $\alpha_{_t}$ must be equal to zero. Otherwise,
$\bar{h}=u_{_t}$ for some $t$ where $[u_{_t}]\in{Irr(S)}$,
which contradicts $(ii)'$.

$(iii)\Rightarrow (i)$. Suppose that
$h$ is any composition of polynomial in $S$. Then $h\in Id(S)$.
Since $Irr(S)$ is a $\mathbf{k}$-basis for $C(B,N)/Id(S)$ and by Lemma \ref{l15},
we have $h\equiv 0\ \ mod(S)$,
i.e.,
 $S$ is a Gr\"{o}bner-Shirshov basis.
\ \ \ \  $\square$

\begin{remark}
 In \cite{BFK04}, a Composition-Diamond lemma for associative conformal algebras
 is established which claim that in  Theorem \ref{t1}, $(i)\Rightarrow(iii)$,
 but not conversely. The reason is that the definitions of
 a Gr\"{o}bner-Shirshov basis in $C(B,N)$ are different, see Remark \ref{remark1}.
\end{remark}

\ \

Given a set $S\subseteq C(B,N)$,  we denote
$\overline{S}=\{\bar{s}\mid s\in S\}$ with respect to
the ordering (\ref{e1.3}).

\begin{lemma}\label{t0}
For any ideal $I$ of $C(B,N)$,
$I$ has  a minimal Gr\"{o}bner-Shirshov basis.
\end{lemma}
\noindent{\bf Proof.}
Define a relation $``\sim"$ on $I$: for any $f,g\in I$,
$$
f\sim g\Leftrightarrow \bar{f}=\bar{g}.
$$
Clearly, $``\sim"$ is an equivalence relation. Let $S_1$ be a set of
representation elements of equivalence classes, say,
$S_1=I/\!\!\sim$. We may assume that $S_1$ is a monic set.  Now we prove $Id(S_1)=I$ and
$S_1$ is a Gr\"{o}bner-Shirshov basis. Clearly, $Id(S_1)\subseteq I$.
Let $0\neq f\in I$. We prove that $f\in Id(S_1)$ by induction on $\bar f$.
If $\bar{f}$ is the minimal element of $\overline{I}$,
then $f\in S_1$.
Otherwise, there exists $g\in S_1$, such that $\bar{f}=\bar{g}$. Now, $f\neq g$,
$\overline{f-g}<\overline{f}$ and $f-g\in I$,
a contradiction.
Suppose that $\bar{f}$ is not the minimal element of $\overline{I}$.
Note that there exists $g\in S_1$ such that $\bar{f}=\bar{g}$. If $f\neq g$, then
$\overline{f-g}<\overline{f}$ and $f-g\in I$. Thus, by induction on $\bar f$, we have $f-g\in Id(S_1)$ and so
$f\in Id(S_1)$. Moreover, by Theorem \ref{t1}, $S_1$ is a Gr\"{o}bner-Shirshov basis for $Id(S_1)=I$.

Note that for any $f,g\in S_1$, if $f\neq g$, then $\bar f\neq \bar g$.

For any $g\in S_1$, define
$$
\triangle_{g}=\{f\in S_1\mid f\neq g,\ \bar{f}=\overline{[u]_{D^ig}}\ \
\mbox{for some normal } S\mbox{-word } [u]_{D^{i}g}\}.
$$
Clearly, if $f\in\triangle_{g}$ and $g\in\triangle_{h}$, then
$\bar{f}>\bar{g}$ and $f\in\triangle_{h}$.
Let
$$
S'=S_1\backslash\cup_{g\in S_1}\triangle_{g}.
$$
Then $Id(S')=Id(S_1)$
and $S'$ is a minimal
Gr\"{o}bner-Shirshov basis for $Id(S')=I$.
In fact, clearly, $Id(S')\subseteq Id(S_1)$.
Suppose that  $f\in S_1$. We prove that $f\in Id(S')$ by induction on $\bar f$.
If $\bar{f}$ is the minimal element of $\overline{S_1}$,
then $f\in S'$. Suppose that  $\bar{f}$ is not the minimal element of $\overline{S_1}$ and $f\notin S'$. Then there exists $g\in S_1$
such that $f\in\triangle_{g}$. So $\bar{g}<\bar{f}$.
By induction on $\bar{f}$, we get that
 there exists $s\in S'$ such that $g\in\triangle_{s}$. So $f\in\triangle_{s}$ and
 $\bar{f}=\overline{[u]_{D^is}}$.
By Theorem \ref{t1} $(ii)'$, we have $f-[u]_{D^is}=\sum\alpha_{_t}[u_{_{t}}]_{_{D^{l_t}s_{_t}}}$,
where each $[u_{_{t}}]_{_{D^{l_t}s_{_t}}}$ is a normal
$S_1$-word and $\overline{[u_{_{t}}]_{_{D^{l_t}s_{_t}}}}\leq\overline{f-[u]_{D^is}}$.
Thus $\overline{s_{_t}}\leq\overline{f-[u]_{D^is}}<\bar{f}$.
By induction on $\bar{f}$, we have $s_{_t}\in Id(S')$. Therefore $f\in Id(S')$.
This shows that $Id(S')=Id(S_1)$.

Let $0\neq f\in Id(S')$. Then $f\in Id(S_1)$. By Theorem \ref{t1},
we have $\bar{f}=\overline{[u]_{D^ig}}$ for some $g\in S_1$,
where $[u]_{D^ig}$ is a normal $S_1$-word.
If $g\in S'$, then $[u]_{D^ig}$ is a normal $S'$-word.
If $g\notin S'$, then there exists $s'\in S'$, such that
$g\in\triangle_{s'},\ \bar{g}=\overline{[w]_{D^ts'}}$.
So $\bar{f}$ is the leading term of some normal $S'$-word $[v]_{D^js'}$.
By Theorem \ref{t1}, $S'$ is a Gr\"{o}bner-Shirshov basis.

This completes the proof.
\ \ \ \  $\square$

\ \

Generally, for an ideal $I$ of $C(B,N)$,
minimal Gr\"{o}bner-Shirshov basis for $I$ is not unique.

\begin{theorem}\label{t*}
For any ideal $I$ of $C(B,N)$,
$I$ has a unique reduced Gr\"{o}bner-Shirshov basis.
\end{theorem}
\noindent{\bf Proof.}
By Lemma \ref{t0}, we may assume that
$I$ has a minimal Gr\"{o}bner-Shirshov basis $S$.

Let $f\in S$. By Lemma \ref{l15}, we have
$$
f-\bar{f}=\sum\alpha_{_k}[u_{_k}]+\sum\alpha_{_j}[u_{_{j}}]_{_{D^{l_j}s_{_j}}},
$$
where  each $[u_{_{j}}]_{_{D^{l_j}s_{_j}}}$ is a normal $S$-word,
$\overline{[u_{_{j}}]_{_{D^{l_j}s_{_j}}}}\leq\overline{f-\bar{f}}<\bar{f}$,
$[u_{_k}]\in Irr(S)$, $u_{_k}\leq\overline{f-\bar{f}}<\bar{f}$,
$\alpha_{_j},\alpha_{_k}\in {\bf k}$ and $s_{_j}\in S$.
Denote $r'_{_{f}}=\sum\alpha_{_k}[u_{_k}]$ and
$$
S'=\{f'\mid f'=\bar{f}+r'_{_{f}},\ f\in S\}.
$$
Then $S'$ is a  reduced Gr\"{o}bner-Shirshov basis for $Id(S)=I$.
In fact,
clearly, $S'\subseteq Id(S)$.
Let $s\in S$. We prove that $s\in Id(S')$ by induction on $\bar s$.
If $\bar{s}$ is the minimal element of $\overline{S}$,
then $s\in S'$. Suppose that  $\bar{s}$ is not the minimal element of $\overline{S}$ and $s\notin S'$.
By Theorem \ref{t1} $(ii)'$, we have
$s-s'=\sum\alpha_{_t}[u_{_{t}}]_{_{D^{l_t}s_{_t}}}$,
where each $[u_{_{t}}]_{_{D^{l_t}s_{_t}}}$ is a normal
$S$-word and
$\overline{[u_{_{t}}]_{_{D^{l_t}s_{_t}}}}\leq\overline{s-s'}$.
Thus $\overline{s_{_t}}\leq\overline{s-s'}<\bar{s}$.
By induction on $\bar{s}$,
we have $s_{_t}\in Id(S')$. So $s\in Id(S')$.
This shows that $Id(S')=Id(S)$.
Since $\overline{S'}=\overline{S}$,\ $Irr(S')=Irr(S)$.
Then $S'$ is
a Gr\"{o}bner-Shirshov basis following from Theorem \ref{t1}.
Clearly, $S'$ is a reduced Gr\"{o}bner-Shirshov basis.

Assume that $S_1$ is another reduced Gr\"{o}bner-Shirshov basis
for ideal $I=Id(S)$. Then $Id(S_1)=Id(S)=Id(S')$.
For any $s_1\in S_1\subseteq Id(S')$, there exists $s'\in S'$
such that $\overline{s_1}=\overline{[u]_{D^is'}}$ by Theorem \ref{t1}.
Similarly, for $s'\in Id(S_1)$, there exists $f_1\in S_1$
such that $\overline{s'}=\overline{[v]_{D^jf_1}}$.
Thus $\overline{s_1}=\overline{[u]_{D^is'}}=\overline{[w]_{D^tf_1}}$
and $\overline{s_1}\geq\overline{ s'}\geq\overline{f_1}$.
Since $S_1$ is a
reduced Gr\"{o}bner-Shirshov basis, we have $\overline{s_1}=\overline{f_1}$.
Thus $\overline{s_1}=\overline{s'}$ and so $\overline{S_1}\subseteq\overline{S'}$.
Similarly, $\overline{S'}\subseteq\overline{S_1}$.
Hence $\overline{S'}=\overline{S_1}$ and so $Irr(S_1)=Irr(S')=Irr(S)$.
Since $S_1$ is a
reduced Gr\"{o}bner-Shirshov basis, for any $s_1\in S_1$, we can assume that $s_1=\overline{s_1}+r_{_{s_1}}$
and $r_{_{s_1}}\in \mathbf{k}Irr(S)$,
where $\mathbf{k}Irr(S)$
is the vector space spanned by $\mathbf{k}$-basis $Irr(S)$. Let $s'\in S'$ with $\overline{s_1}=\overline{s'}$.
Then $s_1-s'=r_{_{s_1}}-r'_{_{s}}\in Id(S)$ and
$r_{_{s_1}}-r'_{_{s}}\in\mathbf{k}Irr(S)$.
By Theorem \ref{t1}, $r_{_{s_1}}=r'_{_{s}}$.
Then $s_1=s'\in S'$. Thus $S_1\subseteq S'$.
Similarly, $S'\subseteq S_1$.
Therefore, $S_1=S'$.
\ \ \ \  $\square$

\ \

\noindent{\bf Shirshov Algorithm} Let $S$ be a monic subset of $C(B,N)$.
If $S$  is not a Gr\"{o}bner-Shirshov basis in $C(B,N)$, then
one can add all nontrivial compositions of polynomials of $S$ to
$S$. Continuing this process repeatedly, we finally obtain a
Gr\"{o}bner-Shirshov basis $S^{c}$ that contains $S$. Such a process
is called Shirshov algorithm.

\section{Applications}

The following proposition gives a characterization of a conformal algebra
$C(B,N|S)$ to have a $\mathbf{k}[D]$-basis under the condition
``$ind(\bar{s})=0$ for any $s\in S$", where  $\mathbf{k}[D]$ is the
polynomial ring over the field  $\mathbf{k}$ with one variable $D$.

\begin{proposition}\label{t2}
Let $S$ be a monic subset of $C(B,N)$ with $ind(\bar{s})=0$ for any $s\in S$.
Then $S$ is a Gr\"{o}bner-Shirshov basis in $C(B,N)$ if and only if
the set
$A=\{[u]\in Irr(S)\mid ind(u)=0\}$ is a $\mathbf{k}[D]$-basis of $C(B,N|S)=C(B,N)/Id(S)$.
\end{proposition}
\noindent{\bf Proof.}
It is clear that $Irr(S)=\bigcup\limits_{[u]\in A}\{[uD^i]\mid i\geq 0\}$.

Suppose that $f\in C(B,N)$. By induction on $\bar f$, we have
\begin{equation*}
f=\sum\limits_{i}\alpha_{_i}[w_{_i}]_{_{D^{l_i}s_i}}
+\sum\limits_{j}\alpha_{_j}D^{l_j}[u_{_j}],
\end{equation*}
where  each $[w_{_i}]_{_{D^{l_i}s_i}}$ is a normal $S$-word,
$\overline{[w_{_i}]_{_{D^{l_i}s_i}}}\leq\bar f,\ [u_{_j}]\in A,\
u_{_j}D^{l_j}\leq \bar f,\ \alpha_{_i},\alpha_{_j}\in \mathbf{k}$
and $s_{_i}\in S$.

Now the result follows from Lemma \ref{l15}.
\ \ \ \  $\square$

\ \

The following corollary follows from Proposition \ref{t2} and Theorem \ref{t1}.

\begin{corollary}
Let $C$ be an associative conformal algebra with a well-ordered
$\mathbf{k}[D]$-basis $\{b_i\}_{i\in I}$, multiplication table
$$
S=\{b_{i(n)}b_j=\Sigma_t\alpha_{ijn}^t b_t,\  \alpha_{ijn}^t\in \mathbf{k}[D],\ i,j\in I,\ 0\leq n<N\}
$$
and a bounded locality
$N=N(b_{i},b_j)$ for any $i,j\in I$.
Then
$
C=C(\{b_i\}_{i\in I}, N| S).
$
\end{corollary}

\ \

Let $S$ be a monic $D$-free subset of $C(B,N)$. Then by Shirshov
Algorithm, we
can get a $D$-free Gr\"{o}bner-Shirshov basis $S^c$ in $C(B,N)$.
Thus by Proposition \ref{t2}, we have the following corollary.

\begin{corollary}(\cite{BFK00})
Let $S$ be a $D$-free subset of monic polynomials in $C(B,N)$.
Then $C(B,N|S)$ is a free $\mathbf{k}[D]$-module.
\end{corollary}

\ \

Now we consider some Lie conformal algebras which are embeddable
into their universal enveloping associative conformal algebras.

\ \

Any associative conformal algebra
$A=(A,(n), n\in\mathbb{Z}_{\geq 0}, D)$
can be made into a Lie conformal algebra
$A^{(-)}=(A,[n],  n\in\mathbb{Z}_{\geq 0}, D)$
by defining new multiplications using conformal commutators:
$$
x_{[n]}y=x_{(n)}y-\{y_{(n)}x\},
$$
where $x,y\in A, n\in\mathbb{Z}_{\geq 0},\
\{y_{(n)}x\}=\sum\limits_{k\geq
0}(-1)^{n+k}\frac{1}{k!}D^{k}(y_{(n+k)}x)$.
The locality function for $A^{(-)}$ is essentially the same as
for $A$, i.e., it is given by
$$
N_{A^{(-)}}(a,b)=max\{N_{A}(a,b), N_{A}(b,a)\}.
$$

Let $L$ be a Lie conformal algebra which is a free
$\mathbf{k}[D]$-module with a  $\mathbf{k}[D]$-basis $B=\{b_i|i\in I\}$ and a
bounded locality $N(b_i,b_j)\leq N$ for all $i,j\in I$.
Let the multiplication table of $L$ in the $\mathbf{k}[D]$-basis $B$ be
$$
b_{i[n]}b_j=\Sigma\alpha_{ijn}^k b_k,\ \alpha_{ijn}^k\in \mathbf{k}[D], \ i,j\in I, n<N.
$$
Then by $\mathcal{U}_{N}(L)$, a universal
enveloping associative conformal algebra of $L$ with respective to
$B$ and $N$, one means the following associative
conformal algebra (see \cite{Ro00}):
$$
\mathcal{U}_{N}(L)=C(B,N| \
b_{i(n)}b_j -\{b_{j(n)}b_i \}-b_{i[n]}b_j=0, \  i,j\in I, \ n<N).
$$

\begin{example}\label{exam1}
Loop Virasoro Lie conformal algebra.
Let
$$
\mathcal{LW}=C_{Lie}(\{L_i\}_{i\in \mathbb{Z}}, N=2\mid
L_{i[0]}L_j=-DL_{i+j},\ L_{i[1]}L_j=-2L_{i+j},\ i,j\in \mathbb{Z})
$$
be the Lie conformal algebra over the complex field $\mathbb{C}$ with a
$\mathbb{C}[D]$-basis $\{L_i\}_{i\in \mathbb{Z}}$. Then
$\mathcal{LW}$ is called loop Virasoro Lie conformal algebra,
see \cite{exLV}.
A universal enveloping associative conformal algebra of
$\mathcal{LW}$ is then given by
$$
\mathcal{U}_{N=2}(\mathcal{LW})=C_{ass}(\{L_i\}_{i\in \mathbb{Z}}, N=2\mid S^{(-)}),
$$
where $S^{(-)}$ consists of
\begin{align*}
 &f_{ij}^0:=L_{i(0)}L_j+L_{j(1)}DL_i-2L_{j(0)}L_i+DL_{i+j},\ \ \ \ i,j\in\mathbb{Z},\\
 &f_{ij}^1:=L_{i(1)}L_j+L_{j(1)}L_i+2L_{i+j},\ \ \ \ i,j\in\mathbb{Z}.
\end{align*}
Then the following statements hold.
\begin{enumerate}
 \item[(i)]\ $C_{ass}(\{L_i\}_{i\in \mathbb{Z}}, N=2\mid S^{(-)})
  =C(\{L_i\}_{i\in \mathbb{Z}}, N=2\mid S_{_1})$,
where $S_{_1}$ consists of the following polynomials:
\begin{align*}
&s_{ij}^0:=L_{i(0)}L_j-L_{0(0)}L_{i+j},\ i\neq 0,\ \ i,j\in\mathbb{Z},\\
&s_{ij}^1:=L_{i(1)}L_j+L_{i+j},\ \ i,j\in\mathbb{Z}.
\end{align*}
 \item[(ii)]\ Define an ordering on the set $\{L_i\}_{i\in \mathbb{Z}}$: for any $i,j\in\mathbb{Z}$,
$$
L_i>L_j\Leftrightarrow |i|>|j| \ \ \mbox{or}\ \ |i|=|j|\ \mbox{and}\ i>j.
$$
Then
$S_{_1}$ is a Gr\"{o}bner-Shirshov basis in
 $C(\{L_i\}_{i\in \mathbb{Z}}, N=2)$.
 \item[(iii)]\ The set
 $$
 Irr(S_1)=\{[\underbrace{L_{_{0}(0)}L_{_{0}(0)}\cdots L_{_{0}(0)}}_{k}D^tL_i]\mid i\in\mathbb{Z},\ t\geq 0,\ k\geq 0\}
 $$
 is a $\mathbb{C}$-basis of
 $\mathcal{U}_{N=2}(\mathcal{LW})$ and the set
 $$
 A=\{[\underbrace{L_{_{0}(0)}L_{_{0}(0)}\cdots L_{_{0}(0)}}_{k}L_i]\mid i\in\mathbb{Z},\  k\geq 0\}
 $$
 is a $\mathbb{C}[D]$-basis of
 $\mathcal{U}_{N=2}(\mathcal{LW})$.
 Moreover, $\mathcal{LW}$ can be embeddable into its universal
enveloping associative conformal algebra $\mathcal{U}_{N=2}(\mathcal{LW})$.
\end{enumerate}
\end{example}
\noindent{\bf Proof.}
(i)\ Since $f_{ii}^{1}=2(L_{i(1)}L_i+L_{2i})$,
we have $s_{ii}^{1}=\frac{1}{2}f_{ii}^{1}\in Id(S^{(-)})$.
Then
\begin{eqnarray*}
s_{0i}^{1}&=&-\frac{1}{2}(f_{i0(1)}^{1}L_0-L_{i(1)}s_{00}^{1}
-L_{0(1)}f_{i0}^{1}-f_{i0}^{1}+s_{00(1)}^{1}L_i)\in Id(S^{(-)}),\\
s_{i0}^{1}&=&f_{i0}^{1}-s_{0i}^{1}\in Id(S^{(-)}),\\
s_{ij}^{1}&=&\frac{1}{2}f_{ji(2)}^{0}L_0+L_{i(1)}f_{j0}^{1}-L_{i(1)}s_{0j}^{1}+f_{i+j,0}^{1}-s_{0,i+j}^{1}\in Id(S^{(-)}),\\
z_{ij}^{0}:&=&L_{i(0)}L_j-L_{j(0)}L_i=f_{ij}^{0}-Ds_{ji}^{1}\in Id(S^{(-)}),\\
s_{ij}^{0}&=&L_{i(1)}z_{j0}^{0}-s_{ij(0)}^{1}L_0-L_{i(0)}s_{j0}^{1}+s_{i0(0)}^{1}L_j+z_{i+j,0}^{0}+L_{i(0)}s_{0j}^{1}\in Id(S^{(-)}).
\end{eqnarray*}
 Therefore $Id(S_{_1})\subseteq Id(S^{(-)})$.
It is obvious that $f_{ij}^0=s_{ij}^0-s_{ji}^0+Ds_{ji}^1$
and $f_{ij}^1=s_{ij}^1+s_{ji}^1$,
so $Id(S^{(-)})\subseteq Id(S_{_1})$. We complete the proof for (i).

(ii)\ All compositions in $S_{_1}$ consist of
\begin{eqnarray*}
&&s_{ij}^0\wedge s_{jk}^0,\ w=L_{i(0)}L_{j(0)}L_k;\ \ \ \
s_{ij}^0\wedge s_{jk}^1,\ w=L_{i(0)}L_{j(1)}L_k;\\
&&s_{ij}^1\wedge s_{jk}^0,\ w=L_{i(1)}L_{j(0)}L_k;\ \ \ \
s_{ij}^1\wedge s_{jk}^1,\ w=L_{i(1)}L_{j(1)}L_k;\\
&&x_{(n)}s_{ij}^0,\ x_{(n)}s_{ij}^1,\ \ x\in\{L_i\}_{i\in\mathbb{Z}},\ n\geq 2.
\end{eqnarray*}
Now we prove that these compositions are trivial modulo $S_{_1}$.

For  $s_{ij}^0\wedge s_{jk}^0,\ w=L_{i(0)}L_{j(0)}L_k$,
\begin{eqnarray*}
(s_{ij}^0, s_{jk}^0)_w
&=&(L_{i(0)}L_j-L_{0(0)}L_{i+j})_{(0)}L_k-
L_{i(0)}(L_{j(0)}L_k-L_{0(0)}L_{j+k})\\
&=&-L_{0(0)}(L_{i+j(0)}L_k)+L_{i(0)}(L_{0(0)}L_{j+k})\\
&=&-L_{0(0)}(L_{i+j(0)}L_k)+(L_{i(0)}L_{0})_{(0)}L_{j+k}\\
&\equiv &-L_{0(0)}(L_{i+j(0)}L_k)+L_{0(0)}(L_{i(0)}L_{j+k})\\
&\equiv &0\ mod(S_{_1}).
\end{eqnarray*}

For $s_{ij}^0\wedge s_{jk}^1,\ w=L_{i(0)}L_{j(1)}L_k$,
we have
\begin{eqnarray*}
(s_{ij}^0, s_{jk}^1)_w&=&(L_{i(0)}L_j-L_{0(0)}L_{i+j})_{(1)}L_k-
L_{i(0)}(L_{j(1)}L_k+L_{j+k})\\
&=&-L_{0(0)}(L_{i+j(1)}L_k)-L_{i(0)}L_{j+k}\\
&\equiv &L_{0(0)}L_{i+j+k}-L_{0(0)}L_{i+j+k}\\
&\equiv &0\ mod(S_{_1}).
\end{eqnarray*}

For $s_{ij}^1\wedge s_{jk}^0,\ w=L_{i(1)}L_{j(0)}L_k$,
we have
\begin{eqnarray*}
(s_{ij}^1, s_{jk}^0)_w
&=&(L_{i(1)}L_j+L_{i+j})_{(0)}L_k-
L_{i(1)}(L_{j(0)}L_k-L_{0(0)}L_{j+k})\\
&=&-L_{i(0)}(L_{j(1)}L_k)+L_{i+j(0)}L_k+(L_{i(1)}L_{0})_{(0)}L_{j+k}+(L_{i(0)}L_{0})_{(1)}L_{j+k}\\
&\equiv &L_{i(0)}L_{j+k}+L_{i+j(0)}L_k-L_{i(0)}L_{j+k}-L_{i(0)}L_{j+k}\\
&\equiv &0\ mod(S_{_1}).
\end{eqnarray*}

For $s_{ij}^1\wedge s_{jk}^1,\ w=L_{i(1)}L_{j(1)}L_k$,
we have
\begin{eqnarray*}
(s_{ij}^1, s_{jk}^1)_w&=&(L_{i(1)}L_j+L_{i+j})_{(1)}L_k-
L_{i(1)}(L_{j(1)}L_k+L_{j+k})\\
&=&L_{i+j(1)}L_k-L_{i(1)}L_{j+k}\\
&\equiv &0\ mod(S_{_1}).
\end{eqnarray*}

For any $t\in\mathbb{Z},\ n\geq 2$, it is clearly that
$L_{t(n)}s_{ij}^1\equiv 0\ mod(S_{_1})$ and
\begin{eqnarray*}
L_{t(n)}s_{ij}^0&=&L_{t(n)}(L_{i(0)}L_j-L_{0(0)}L_{i+j})\\
&=&nL_{t(n-1)}(L_{i(1)}L_j-L_{0(1)}L_{i+j})\\
&\equiv &nL_{t(n-1)}(L_{i+j}-L_{i+j})\\
&\equiv &0\ mod(S_{_1}).
\end{eqnarray*}

Therefore, $S_{_1}$ is a Gr\"{o}bner-Shirshov basis.

(iii)\ By Theorem \ref{t1} and  Proposition \ref{t2},
we can get the result.\ \ \ \  $\square$

\begin{example}\label{exam2}
Loop Heisenberg-Virasoro Lie conformal algebra.
Let
\begin{eqnarray*}
\mathcal{CHV}&=&C_{Lie}(\{L_i,H_j\}_{i,j\in \mathbb{Z}},\  N=2\mid
L_{i[0]}L_j=DL_{i+j},\ L_{i[1]}L_j=2L_{i+j},\\
&&\ \ \ \ \ \ \  L_{i[0]}H_j=DH_{i+j}, L_{i[0]}H_j=DH_{i+j},\ L_{i[1]}H_j=H_{i+j},\\
&&\ \ \ \ \ \ \ H_{i[0]}H_j=H_{i[1]}H_j=0,\ i,j\in \mathbb{Z})
\end{eqnarray*}
be the Lie conformal algebra over the complex field $\mathbb{C}$
with a  $\mathbb{C}[D]$-basis $\{L_i, H_j\}_{i, j\in \mathbb{Z}}$. Then
$\mathcal{CHV}$ is called loop Heisenberg-Virasoro Lie conformal algebra,
see \cite{exLH}.
A universal enveloping associative conformal algebra of
$\mathcal{CHV}$ is then given by
$$
\mathcal{U}_{N=2}(\mathcal{CHV})=C_{ass}(\{L_i, H_j\}_{i,j\in\mathbb{Z}},\ N=2\mid S^{(-)}),
$$
where $S^{(-)}$ consists of
\begin{align*}
&f_{ij}^0:=L_{i(0)}L_j+L_{j(1)}DL_i-2L_{j(0)}L_i-DL_{i+j}, \\
&f_{ij}^1:=L_{i(1)}L_j+L_{j(1)}L_i-2L_{i+j}, \\
&g_{ij}^0:=L_{i(0)}H_j+H_{j(1)}DL_i-2H_{j(0)}L_i-DH_{i+j},  \\
&g_{ij}^1:=L_{i(1)}H_j+H_{j(1)}L_i-H_{i+j},\\
&p_{ij}^0:=L_{j(1)}DH_i-2L_{j(0)}H_i+H_{i(0)}L_j, \\
&h_{ij}^0:=H_{i(0)}H_j+H_{j(1)}DH_i-2H_{j(0)}H_i, \\
&h_{ij}^1:=H_{i(1)}H_j+H_{j(1)}H_i
\end{align*}
for any $i,j\in \mathbb{Z}$.  Then the following statements hold.
\begin{enumerate}
 \item[(i)]\ $C_{ass}(\{L_i, H_j\}_{i,j\in \mathbb{Z}},\ N=2\mid S^{(-)})
 =C(\{L_i, H_j\}_{i,j\in \mathbb{Z}},\ N=2\mid S_{_1})$,
where $S_{_1}$ consists of the following polynomials:
\begin{align*}
&s_{ij}^0:=L_{i(0)}L_j-L_{0(0)}L_{i+j},\ i\neq 0,\\
&s_{ij}^1:=L_{i(1)}L_j-L_{i+j}, \\
&g_{ij}^0=L_{i(0)}H_j+H_{j(1)}DL_i-2H_{j(0)}L_i-DH_{i+j},  \\
&g_{ij}^1=L_{i(1)}H_j+H_{j(1)}L_i-H_{i+j},\\
&q_{ijk}^0:=H_{i(0)}L_{j+k}-H_{i+j(0)}L_{k}+H_{j(0)}L_{i+k}-H_{0(0)}L_{i+j+k},\\
& \ \ \ \ \ \ \ \ \ \ \ \ \ \ \ \ \ \ \ \ \ \ \ \
|i|\geq |j|,\ i>0>j\  \mbox{ or }\ i>j>0\  \mbox{ or }\ i<j<0,\\
&q_{ij}^1:=H_{i(1)}L_j-H_{0(1)}L_{i+j},\ i\neq 0,\\
&r_{ij}^0:=H_{i(0)}H_j-H_{0(0)}H_{i+j},\ i\neq 0,\\
&r_{ij}^1:=H_{i(1)}H_j,
\end{align*}
where $i,j,k\in\mathbb{Z}$.
 \item[(ii)]\ Define an ordering set $\{L_i, H_j\}_{i, j\in \mathbb{Z}}$:
 for any $i,j, k, t\in \mathbb{Z}$, $L_k>H_t $ and
$$
L_i>L_j\ (H_i>H_j)\Leftrightarrow |i|>|j| \ \ \mbox{or}\ \ |i|=|j|\ \mbox{and}\ i>j.
$$

Then
 $S_{_1}$ is a Gr\"{o}bner-Shirshov basis in
 $C(\{L_i, H_j\}_{i,j\in \mathbb{Z}},\ N=2)$.
 \item[(iii)]\ The set
 \begin{eqnarray*}
 Irr(S_1)&=&\{[\underbrace{H_{_{0}(0)}H_{_{0}(0)}\cdots H_{_{0}(n)}}_{k}
 \underbrace{L_{_{0}(0)}L_{_{0}(0)}\cdots L_{_{0}(0)}}_{l}D^tL_i],\ H_{_{0}(0)}D^tH_i,\\
 && H_{_{-1}(0)}D^tL_i,\ D^tH_i,  \ i\in\mathbb{Z},\ t\geq 0,k\geq 0,l\geq 0,\ 0\leq n\leq 1\}
 \end{eqnarray*}
 is a $\mathbb{C}$-basis of
 $\mathcal{U}_{N=2}(\mathcal{CHV})$ and the set
 \begin{eqnarray*}
 A&=&\{[\underbrace{H_{_{0}(0)}H_{_{0}(0)}\cdots H_{_{0}(n)}}_{k}
 \underbrace{L_{_{0}(0)}L_{_{0}(0)}\cdots L_{_{0}(0)}}_{l}L_i],\  H_{_{0}(0)}H_i,\ H_i,\\
 && H_{_{-1}(0)}L_i,\ \ i\in\mathbb{Z},\ k\geq 0,\ l\geq 0,\ 0\leq n\leq 1\}
 \end{eqnarray*}
 is a $\mathbb{C}[D]$-basis of
 $\mathcal{U}_{N=2}(\mathcal{CHV})$.
 Moreover, $\mathcal{CHV}$ can be embeddable into its universal
enveloping associative conformal algebra $\mathcal{U}_{N=2}(\mathcal{CHV})$.
\end{enumerate}
\end{example}
\noindent{\bf Proof.}
(i)\ By calculation, we have
\begin{eqnarray*}
s_{ii}^{1}&=&\frac{1}{2}f_{ii}^{1}\in Id(S^{(-)}),\\
r_{ii}^{1}&=&\frac{1}{2}h_{ii}^{1}\in Id(S^{(-)}),\\
s_{i0}^{1}&=&-\frac{1}{2}(f_{i0(1)}^{1}L_0-L_{i(1)}s_{00}^{1}-
L_{0(1)}f_{i0}^{1}-f_{i0}^{1}+s_{00(1)}^{1}L_i)\in Id(S^{(-)}),\\
s_{0i}^{1}&=&f_{i0}^{1}-s_{i0}^{1}\in Id(S^{(-)}),\\
s_{ij}^{1}&=&-\frac{1}{2}f_{ji(2)}^{0}L_0-
L_{i(1)}s_{j0}^{1}+s_{i+j,0}^{1}\in Id(S^{(-)}),\\
z_{ij}^{0}:&=&L_{i(0)}L_j-L_{j(0)}L_i=f_{ij}^{0}-
Ds_{ji}^{1}\in Id(S^{(-)}),\\
s_{ij}^{0}&=&-L_{i(1)}z_{j0}^{0}+s_{ij(0)}^{1}L_0+
L_{i(0)}s_{j0}^{1}-s_{i0(0)}^{1}L_j+z_{i+j,0}^{0}-
L_{i(0)}s_{0j}^{1}\in Id(S^{(-)}),\\
r_{ij}^{1}&=&\frac{1}{2}(u_{ij}^{1}
-u_{ji}^{1}+h_{ji}^{1})\in Id(S^{(-)}),\\
v_{ij}^{0}:&=&H_{i(0)}H_j-H_{j(0)}H_i
=h_{ij}^{0}-Dr_{ji}^{1}\in Id(S^{(-)}),\\
r_{ij}^{0}&=&-g_{i0(1)}^{1}H_j+L_{i(1)}v_{0j}^{0}
-L_{i(0)}r_{0j}^{1}+H_{0(1)}g_{ij}^{0}-H_{0(0)}g_{ij}^{1}+g_{ij(0)}^{1}H_0 \\
&&+L_{i(0)}r_{j0}^{1}-r_{0j(1)}^1DL_i
+2r_{0j(0)}^1L_i+Dr_{0,i+j}^1+r_{j0(1)}^1DL_i-2r_{j0(0)}^1L_i\\
&&+2v_{j0(1)}^{0}DL_i+3v_{0j(1)}^{0}L_i
-Dr_{ji}^1+v_{i+j,0}^0 \in Id(S^{(-)}),\\
q_{ij}^{1}&=&\frac{1}{2}L_{j(2)}g_{i0}^0
-s_{ji}^{1}-g_{i+j,0}^{1}+g_{ji}^{1}\in Id(S^{(-)}),\\
q_{ijk}^{0}&=&-g_{ij(1)}^{0}L_k+L_{i(0)}q_{jk}^{1}-q_{ji(0)}^{1}L_k
-H_{0(1)}s_{i+j,k}^{0}+H_{0(1)}s_{i,j+k}^{0}\in Id(S^{(-)}).
\end{eqnarray*}
Therefore $Id(S_1)\subseteq Id(S^{(-)})$.
It is obvious that
\begin{eqnarray*}
f_{ij}^0&=&s_{ij}^0-s_{ji}^0+Ds_{ji}^1, \ \ f_{ij}^1=s_{ij}^1+s_{ji}^1,\ \  h_{ij}^0=r_{ij}^0-r_{ji}^0+Dr_{ji}^1, \\
h_{ij}^1&=&r_{ij}^1+r_{ji}^1,\ \  p_{ij}^0=Dg_{ji}^1-g_{ji}^0.
\end{eqnarray*}
Therefore $Id(S^{(-)})\subseteq Id(S_{_1})$.
Hence $Id(S^{(-)})=Id(S_{_1})$.

(ii)\ Note that if $|i|\geq |j|$ and $i>0>j$, then $\overline{q_{ijk}^0}=H_{i(0)}L_{j+k}$;
if $i>j>0$ or $i<j<0$, then $\overline{q_{ijk}^0}=H_{i+j(0)}L_{k}$.

All compositions in $S_{_1}$ are related to
\begin{eqnarray*}
&&s_{ij}^0\wedge s_{jk}^0,\ w=L_{i(0)}L_{j(0)}L_k; \ \ \ \
s_{ij}^0\wedge s_{jk}^1,\ w=L_{i(0)}L_{j(1)}L_k;\\
&&s_{ij}^0\wedge g_{jk}^0,\ w=L_{i(0)}L_{j(0)}H_k;\ \ \ \
s_{ij}^0\wedge g_{jk}^1,\ w=L_{i(0)}L_{j(1)}H_k;\\
&&s_{ij}^1\wedge s_{jk}^0,\ w=L_{i(1)}L_{j(0)}L_k;\ \ \ \
s_{ij}^1\wedge s_{jk}^1,\ w=L_{i(1)}L_{j(1)}L_k;\\
&&s_{ij}^1\wedge g_{jk}^0,\ w=L_{i(1)}L_{j(0)}H_k;\ \ \ \
s_{ij}^1\wedge g_{jk}^1,\ w=L_{i(1)}L_{j(1)}H_k;\\
&&g_{ij}^0\wedge q_{jkm}^0,\ w=L_{i(0)}H_{j(0)}L_{k+m};\ \ \ \
g_{i,j+k}^0\wedge q_{jkm}^0,\ w=L_{i(0)}H_{j+k(0)}L_{m};\\
&&g_{ij}^0\wedge q_{jk}^1,\ w=L_{i(0)}H_{j(1)}L_k;\ \ \ \
g_{ij}^0\wedge r_{jk}^0,\ w=L_{i(0)}H_{j(0)}H_k;\\
&&g_{ij}^0\wedge r_{jk}^1,\ w=L_{i(0)}H_{j(1)}H_k;\ \ \ \
g_{ij}^1\wedge q_{jkm}^0,\ w=L_{i(1)}H_{j(0)}L_{k+m};\\
&&g_{i,j+k}^1\wedge q_{jkm}^0,\ w=L_{i(1)}H_{j+k(0)}L_m;\ \ \ \
g_{ij}^1\wedge q_{jk}^1,\ w=L_{i(1)}H_{j(1)}L_k;\\
&&g_{ij}^1\wedge r_{jk}^0,\ w=L_{i(1)}H_{j(0)}H_k;\ \ \ \
g_{ij}^1\wedge r_{jk}^1,\ w=L_{i(1)}H_{j(1)}H_k;\\
&&q_{ijk}^0\wedge s_{j+k,m}^0,\ w=H_{i(0)}L_{j+k(0)}L_m;\ \ \ \
q_{ijk}^0\wedge s_{km}^0,\ w=H_{i+j(0)}L_{k(0)}L_m;\\
&&q_{ijk}^0\wedge s_{j+k,m}^1,\ w=H_{i(0)}L_{j+k(1)}L_m;\ \ \ \
q_{ijk}^0\wedge s_{km}^1,\ w=H_{i+j(0)}L_{k(1)}L_m;\\
&&q_{ijk}^0\wedge g_{j+k,m}^0,\ w=H_{i(0)}L_{j+k(0)}H_m;\ \ \ \
q_{ijk}^0\wedge g_{km}^0,\ w=H_{i+j(0)}L_{k(0)}H_m;\\
&&q_{ijk}^0\wedge g_{j+k,m}^1,\ w=H_{i(0)}L_{j+k(1)}H_m;\ \ \ \
q_{ijk}^0\wedge g_{km}^1,\ w=H_{i+j(0)}L_{k(1)}H_m;\\
&&q_{ij}^1\wedge s_{jk}^0,\ w=H_{i(1)}L_{j(0)}L_k;\ \ \ \
q_{ij}^1\wedge s_{jk}^1,\ w=H_{i(1)}L_{j(1)}L_k;\\
&&q_{ij}^1\wedge g_{jk}^0,\ w=H_{i(1)}L_{j(0)}H_k;\ \ \ \
q_{ij}^1\wedge g_{jk}^1,\ w=H_{i(1)}L_{j(1)}H_k;\\
&&r_{ij}^0\wedge q_{jkm}^0,\ w=H_{i(0)}H_{j(0)}L_{k+m};\ \ \ \
r_{i,j+k}^0\wedge q_{jkm}^0,\ w=H_{i(0)}H_{j+k(0)}L_{m};\\
&&r_{ij}^0\wedge q_{jk}^1,\ w=H_{i(0)}H_{j(1)}L_k;\ \ \ \
r_{ij}^0\wedge r_{jk}^0,\ w=H_{i(0)}H_{j(0)}H_k;\\
&&r_{ij}^0\wedge r_{jk}^1,\ w=H_{i(0)}H_{j(1)}H_k;\ \ \ \
r_{ij}^1\wedge q_{jkm}^0,\ w=H_{i(1)}H_{j(0)}L_{k+m};\\
&&r_{i,j+k}^1\wedge q_{jkm}^0,\ w=H_{i(1)}H_{j+k(0)}L_{m};\ \ \ \
r_{ij}^1\wedge q_{jk}^1,\ w=H_{i(1)}H_{j(1)}L_k;\\
&&r_{ij}^1\wedge r_{jk}^0,\ w=H_{i(1)}H_{j(0)}H_k;\ \ \ \
r_{ij}^1\wedge r_{jk}^1,\ w=H_{i(1)}H_{j(1)}H_k;\\
&&x_{(n)}s;\ \ \ \
g_{ij(n)}^0x; \ \ \ \mbox{where }\ s\in S_1,\ x\in\{L_i, H_j\}_{i,j\in\mathbb{Z}},\ n\geq 2.
\end{eqnarray*}

Now we  prove that these compositions are trivial modulo $S_{_1}$.

It is similar to Example \ref{exam1} that for
$s_{ij}^0\wedge s_{jk}^0, w=L_{i(0)}L_{j(0)}L_k$;
$s_{ij}^0\wedge s_{jk}^1, w=L_{i(0)}L_{j(1)}L_k;\
s_{ij}^1\wedge s_{jk}^0, w=L_{i(1)}L_{j(0)}L_k;\
s_{ij}^1\wedge s_{jk}^1, w=L_{i(1)}L_{j(1)}L_k$;
$r_{ij}^0\wedge r_{jk}^0, w=H_{i(0)}H_{j(0)}H_k$;
$r_{ij}^0\wedge r_{jk}^1, w=H_{i(0)}H_{j(1)}H_k;\
r_{ij}^1\wedge s_{jk}^0, w=H_{i(1)}H_{j(0)}H_k$;
$r_{ij}^1\wedge r_{jk}^1, w=H_{i(1)}H_{j(1)}H_k$;
the corresponding compositions are trivial modulo $S_{_1}$.

For $s_{ij}^0\wedge g_{jk}^1,\ w=L_{i(0)}L_{j(1)}H_k$, we have $i\neq 0$ and \begin{eqnarray*}
&&(s_{ij}^0, g_{jk}^1)_w\\
&=&(L_{i(0)}L_j-L_{0(0)}L_{i+j})_{(1)}H_k
-L_{i(0)}(L_{j(1)}H_k+H_{k(1)}L_j-H_{j+k})\\
&=&-L_{0(0)}(L_{i+j(1)}H_k)-(L_{i(0)}H_{k})_{(1)}L_j+L_{i(0)}H_{j+k}\\
&\equiv & L_{0(0)}(H_{k(1)}L_{i+j}-H_{i+j+k})+(H_{k(1)}DL_i
-2H_{k(0)}L_i-DH_{i+k})_{(1)}L_j\\
&&-H_{j+k(1)}DL_i+2H_{j+k(0)}L_i+DH_{i+j+k}\\
&\equiv &-(H_{k(1)}DL_0-2H_{k(0)}L_0-DH_{k})_{(1)}L_{i+j}
+H_{i+j+k(1)}DL_0-2H_{i+j+k(0)}L_0\\
&&-DH_{i+j+k}-H_{k(1)}(L_{i(0)}L_j)+2H_{k(0)}(L_{i(1)}L_j)
-2H_{k(0)}(L_{i(1)}L_j)\\
&&+H_{i+k(0)}L_j-H_{j+k(1)}DL_i+2H_{j+k(0)}L_i+DH_{i+j+k}\\
&\equiv &H_{k(1)}(L_{0(0)}L_{i+j})-2H_{k(0)}(L_{0(1)}L_{i+j})
+2H_{k(0)}(L_{0(1)}L_{i+j})-H_{k(0)}L_{i+j}\\
&&+H_{i+j+k(1)}DL_0-2H_{i+j+k(0)}L_0-DH_{i+j+k}
-H_{k(1)}(L_{0(0)}L_{i+j})\\
&&+H_{i+k(0)}L_j-H_{j+k(1)}DL_i+2H_{j+k(0)}L_i+DH_{i+j+k}\\
&\equiv &-H_{k(0)}L_{i+j}+H_{i+j+k(1)}DL_0
-2H_{i+j+k(0)}L_0+H_{i+k(0)}L_j-H_{j+k(1)}DL_i\\
&&+2H_{j+k(0)}L_i\\
&\equiv &D(H_{i+j+k(1)}L_0-H_{j+k(1)}L_i)
-H_{k(0)}L_{i+j}-H_{i+j+k(0)}L_0+H_{i+k(0)}L_j\\
&&+H_{j+k(0)}L_i\\
&\equiv &H_{i+k(0)}L_j-H_{k(0)}L_{i+j}+H_{j+k(0)}L_i
-H_{i+j+k(0)}L_0\\
&\equiv & \left\{
\begin{aligned}
&H_{i+k(0)}L_0-H_{k(0)}L_{i}+H_{k(0)}L_i
-H_{i+k(0)}L_0,\ \  j=0   \\
&H_{0(0)}L_j-H_{k(0)}L_{j+i}+H_{j+k(0)}L_i
-H_{j(0)}L_0,\ \   j\neq 0,\ i+k=0\\
&-H_{j(0)}L_{i+k}+H_{0(0)}L_{i+j+k}
-H_{k(0)}L_{i+j}+H_{j+k(0)}L_i,\ \   j\neq 0,\ i+k\neq 0
\end{aligned}
\right. \\
&\equiv &0\ mod(S_{_1}).
\end{eqnarray*}

For $s_{ij}^0\wedge g_{jk}^0,\ w=L_{i(0)}L_{j(0)}H_k$, we have
\begin{eqnarray*}
&&(s_{ij}^0, g_{jk}^0)_w\\
&=&(L_{i(0)}L_j-L_{0(0)}L_{i+j})_{(0)}H_k
-L_{i(0)}(L_{j(0)}H_k+H_{k(1)}DL_j-2H_{k(0)}L_j\\
&&-DH_{j+k})\\
&=&-L_{0(0)}(L_{i+j(0)}H_k)-(L_{i(0)}H_{k})_{(1)}DL_{j}+
2(L_{i(0)}H_{k})_{(0)}L_{j}+L_{i(0)}DH_{j+k}\\
&\equiv &L_{0(0)}(H_{k(1)}DL_{i+j}-2H_{k(0)}L_{i+j}-DH_{i+j+k})+
(H_{k(1)}DL_i-2H_{k(0)}L_i\\
&&-DH_{i+k})_{(1)}DL_{j}-2(H_{k(1)}DL_i-2H_{k(0)}L_i-
DH_{i+k})_{(0)}L_{j}\\
&&-D(H_{j+k(1)}DL_i-2H_{j+k(0)}L_i-DH_{i+j+k})\\
\end{eqnarray*}

\begin{eqnarray*}
&\equiv & (L_{0(0)}H_{k})_{(1)}DL_{i+j}-2(L_{0(0)}H_{k})_{(0)}L_{i+j}-
L_{0(0)}DH_{i+j+k}\\
&&-H_{k(1)}(L_{i(0)}DL_{j})+H_{i+k(0)}DL_{j}-2H_{k(0)}(L_{i(0)}L_{j})+4H_{k(0)}(L_{i(0)}L_{j}) \\
&&-H_{j+k(1)}D^2L_i+H_{j+k(0)}DL_i
+2H_{j+k(0)}DL_i+D^2H_{i+j+k} \\
&\equiv &-(H_{k(1)}DL_0-2H_{k(0)}L_0-DH_{k})_{(1)}DL_{i+j} +2(H_{k(1)}DL_0-2H_{k(0)}L_0 \\
&&-DH_{k})_{(0)}L_{i+j}+D(H_{i+j+k(1)}DL_0-
2H_{i+j+k(0)}L_0-DH_{i+j+k})\\
&&-H_{k(1)}D(L_{0(0)}L_{i+j})+
H_{i+k(0)}DL_{j}+2H_{k(0)}(L_{0(0)}L_{i+j})-H_{j+k(1)}D^2L_i\\
&&+H_{j+k(0)}DL_i+2H_{j+k(0)}DL_i+D^2H_{i+j+k} \\
&\equiv &H_{k(1)}(L_{0(0)}DL_{i+j})-2H_{k(0)}(L_{0(1)}DL_{i+j})+
2H_{k(0)}(L_{0(1)}DL_{i+j}) \\
&&-H_{k(0)}DL_{i+j}+2H_{k(0)}(L_{0(0)}L_{i+j})-4H_{k(0)}(L_{0(0)}L_{i+j})+H_{i+j+k(1)}D^2L_0\\
&&-H_{i+j+k(0)}DL_0-2H_{i+j+k(0)}DL_0-H_{k(1)}(L_{p(0)}DL_{i+j})+
H_{i+k(0)}DL_{j}\\
&&+2H_{k(0)}(L_{0(0)}L_{i+j})-H_{j+k(1)}D^2L_i+H_{j+k(0)}DL_i+2H_{j+k(0)}DL_i\\
&\equiv &-H_{k(0)}DL_{i+j}+H_{i+j+k(1)}D^2L_0-H_{i+j+k(0)}DL_0
-2H_{i+j+k(0)}DL_0\\
&&+H_{i+k(0)}DL_{j}-H_{j+k(1)}D^2L_i+H_{j+k(0)}DL_i+2H_{j+k(0)}DL_i\\
&\equiv &D(H_{i+k(0)}L_{j}-H_{k(0)}L_{i+j}+H_{j+k(0)}L_i
-H_{i+j+k(0)}L_0)+D^2(H_{i+j+k(1)}L_0)\\
&&-D^2(H_{j+k(1)}L_i)\\
&\equiv &0\ mod(S_{_1}).
\end{eqnarray*}

For $s_{ij}^1\wedge g_{jk}^0,\ w=L_{i(1)}L_{j(0)}H_k$, we have
\begin{eqnarray*}
&&(s_{ij}^1, g_{jk}^0)_w\\
&=&(L_{i(1)}L_j-L_{i+j})_{(0)}H_k-L_{i(1)}(L_{j(0)}H_k
+H_{k(1)}DL_j-2H_{k(0)}L_j-DH_{j+k})\\
&=&-L_{i(0)}(L_{j(1)}H_k)-L_{i+j(0)}H_k-(L_{i(1)}H_k)_{(1)}DL_j
-(L_{i(0)}H_k)_{(2)}DL_j\\
&&+2(L_{i(1)}H_k)_{(0)}L_j+2(L_{i(0)}H_k)_{(1)}L_j+D(L_{i(1)}H_{j+k})+L_{i(0)}H_{j+k}\\
&\equiv &L_{i(0)}(H_{k(1)}L_j-H_{j+k})+H_{k(1)}DL_{i+j}
-2H_{k(0)}L_{i+j}-DH_{i+j+k}\\
&&+(H_{k(1)}L_i)_{(1)}DL_j-H_{i+k(1)}DL_j+(H_{k(1)}DL_i-2H_{k(0)}L_i-DH_{i+k})_{(2)}DL_j\\
&&-2(H_{k(1)}L_i-H_{i+k})_{(0)}L_j-2(H_{k(1)}DL_i-2H_{k(0)}L_i-DH_{i+k})_{(1)}L_j\\
&&-D(H_{j+k(1)}L_i-H_{i+j+k})+L_{i(0)}H_{j+k}\\
&\equiv &(L_{i(0)}H_{k})_{(1)}L_j+D(H_{k(1)}L_{i+j})-H_{k(0)}L_{i+j}
+H_{k(1)}(L_{i(1)}DL_j)\\
&&-D(H_{i+k(1)}L_j)-H_{k(0)}(L_{i(2)}DL_j)+H_{i+k(0)}L_j-2H_{k(1)}(L_{i(1)}DL_j)\\
&&+2H_{i+k(1)}DL_j+H_{k(0)}(L_{i(2)}DL_j)-2H_{k(1)}(L_{i(0)}L_j)
+2H_{k(1)}(L_{i(0)}L_j)\\
&&-2H_{k(0)}(L_{i(1)}L_j)+4H_{k(0)}(L_{i(1)}L_j)-2H_{i+k(0)}L_j-D(H_{j+k(1)}L_i)\\
\end{eqnarray*}

\begin{eqnarray*}
&\equiv &-(H_{k(1)}DL_i-2H_{k(0)}L_i-DH_{i+k})_{(1)}L_j
+D(H_{k(1)}L_{i+j})-H_{k(1)}DL_{i+j}\\
&&-H_{k(1)}(L_{i(0)}L_j)-D(H_{i+k(1)}L_j)-H_{i+k(0)}L_j+2H_{i+k(1)}DL_j+H_{k(0)}L_{i+j}\\
&&-D(H_{j+k(1)}L_i)\\
&\equiv &D(H_{i+k(1)}L_{j}-H_{j+k(1)}L_{i})\\
&\equiv &0\ mod(S_{_1}).
\end{eqnarray*}

For $s_{ij}^1\wedge g_{jk}^1,\ w=L_{i(1)}L_{j(1)}H_k$, we have
\begin{eqnarray*}
&&(s_{ij}^1, g_{jk}^1)_w\\
&=&(L_{i(1)}L_j-L_{i+j})_{(1)}H_k-
L_{i(1)}(L_{j(1)}H_k+H_{k(1)}L_j-H_{j+k})\\
&=&-L_{i+j(1)}H_k-(L_{i(1)}H_k)_{(1)}L_j+L_{i(1)}H_{j+k}\\
&\equiv &H_{k(1)}L_{i+j}-H_{i+j+k}+(H_{k(1)}L_{i}
-H_{i+k})_{(1)}L_j-H_{j+k(1)}L_{i}+H_{i+j+k}\\
&\equiv &(H_{k(1)}L_{i+j}-H_{i+k(1)}L_{j})+(H_{k(1)}L_{i+j}-H_{j+k(1)}L_{i})\\
&\equiv &0\ mod(S_{_1}).
\end{eqnarray*}

For $g_{ij}^0\wedge q_{jkm}^0,\ w=L_{i(0)}H_{j(0)}L_{k+m}$, we have
\begin{eqnarray*}
&&(g_{ij}^0, q_{jkm}^0)_w\\
&=&(L_{i(0)}H_j+H_{j(1)}DL_i-2H_{j(0)}L_i-DH_{i+j})_{(0)}L_{k+m}-
L_{i(0)}(H_{j(0)}L_{k+m}\\
&&+H_{k(0)}L_{j+m}-H_{0(0)}L_{j+k+m}-H_{j+k(0)}L_{m})\\
&=&H_{j(0)}(L_{i(0)}L_{k+m})-2H_{j(0)}(L_{i(0)}L_{k+m})
-(L_{i(0)}H_{k})_{(0)}L_{j+m}\\
&&+(L_{i(0)}H_{0})_{(0)}L_{j+k+m}+(L_{i(0)}H_{j+k})_{(0)}L_{m}\\
&\equiv &-H_{j(0)}(L_{i(0)}L_{k+m})-H_{k(0)}(L_{i(0)}L_{j+m})+H_{0(0)}(L_{i(0)}L_{j+k+m})\\
&&+H_{j+k(0)}(L_{i(0)}L_{m})\\
&\equiv &(-H_{j(0)}L_{i+k}-H_{k(0)}L_{i+j}+H_{0(0)}L_{i+j+k}+H_{j+k(0)}L_{i})_{(0)}L_{m}\\
&\equiv &0\ mod(S_{_1}).
\end{eqnarray*}

For $g_{i,j+k}^0\wedge q_{jkm}^0,\ w=L_{i(0)}H_{j+k(0)}L_{m}$, we have
\begin{eqnarray*}
&&(g_{i,j+k}^0, -q_{jkm}^0)_w\\
&=&(L_{i(0)}H_{j+k}+H_{j+k(1)}DL_i-2H_{j+k(0)}L_i-DH_{i+j+k})_{(0)}L_{m}\\
&&-
L_{i(0)}(H_{j+k(0)}L_{m}-H_{j(0)}L_{k+m}-H_{k(0)}L_{j+m}+H_{0(0)}L_{j+k+m})\\
\end{eqnarray*}

\begin{eqnarray*}
&=&H_{j+k(0)}(L_{i(0)}L_{m})-2H_{j+k(0)}(L_{i(0)}L_{m})
+(L_{i(0)}H_{k})_{(0)}L_{j+m}\\
&&-(L_{i(0)}H_{0})_{(0)}L_{j+k+m}+(L_{i(0)}H_{j})_{(0)}L_{k+m}\\
&\equiv &-H_{j+k(0)}(L_{i(0)}L_{m})+H_{k(0)}(L_{i(0)}L_{j+m})-H_{0(0)}(L_{i(0)}L_{j+k+m})\\
&&+H_{j(0)}(L_{i(0)}L_{k+m})\\
&\equiv &(-H_{j+k(0)}L_{0}+H_{k(0)}L_{j}-H_{0(0)}L_{j+k}+H_{j(0)}L_{k})_{(0)}L_{i+m}\\
&\equiv &0\ mod(S_{_1}).
\end{eqnarray*}

For $g_{ij}^0\wedge q_{jk}^1,\ w=L_{i(0)}H_{j(1)}L_k$, we have
\begin{eqnarray*}
&&(g_{ij}^0, q_{jk}^1)_w\\
&=&(L_{i(0)}H_j+H_{j(1)}DL_i-2H_{j(0)}L_i-DH_{i+j})_{(1)}L_k
-L_{i(0)}(H_{j(1)}L_k-H_{0(1)}L_{j+k})\\
&=&-H_{j(1)}(L_{i(0)}L_k)+2H_{j(0)}(L_{i(1)}L_k)-2H_{j(0)}(L_{i(1)}L_k)
+H_{i+j(0)}L_k\\
&&+(L_{i(0)}H_{0})_{(1)}L_{j+k}\\
&\equiv &-(H_{j(1)}L_{i})_{(0)}L_k-H_{j(0)}L_{i+k}+H_{i+j(0)}L_k
-(H_{0(1)}DL_i-2H_{0(0)}L_i\\
&&-DH_{i})_{(1)}L_{j+k}\\
&\equiv &-(H_{0(1)}L_{i+j})_{(0)}L_k-H_{j(0)}L_{i+k}+H_{i+j(0)}L_k
+H_{0(1)}(L_{i(0)}L_{j+k})-H_{i(0)}L_{j+k}\\
&\equiv &H_{0(0)}L_{i+j+k}-H_{j(0)}L_{i+k}+H_{i+j(0)}L_k
-H_{i(0)}L_{j+k}\\
&\equiv &0\ mod(S_{_1}).
\end{eqnarray*}

For $g_{ij}^0\wedge r_{jk}^0,\ w=L_{i(0)}H_{j(0)}H_k$, we have
\begin{eqnarray*}
&&(g_{ij}^0, r_{jk}^0)_w\\
&=&(L_{i(0)}H_j+H_{j(1)}DL_i-2H_{j(0)}L_i-DH_{i+j})_{(0)}H_k
-L_{i(0)}(H_{j(0)}H_k-H_{0(0)}H_{j+k})\\
&=&-H_{j(0)}(L_{i(0)}H_k)+(L_{i(0)}H_{0})_{(0)}H_{j+k}\\
&\equiv &H_{j(0)}(H_{k(1)}DL_i-2H_{k(0)}L_i-DH_{i+k})
-(H_{0(1)}DL_i-2H_{0(0)}L_i-DH_{i})_{(0)}H_{j+k}\\
&\equiv &(H_{j(0)}H_{k})_{(1)}DL_i-2(H_{j(0)}H_{k})_{(0)}L_i
-H_{j(0)}DH_{i+k}+H_{0(0)}(L_{i(0)}H_{j+k})\\
&\equiv &(H_{j(0)}H_{k})_{(1)}DL_i-2(H_{j(0)}H_{k})_{(0)}L_i-H_{j(0)}DH_{i+k}\\
&&-H_{0(0)}(H_{j+k(1)}DL_i-2H_{j+k(0)}L_i-DH_{i+j+k})\\
&\equiv &0\ mod(S_{_1}).
\end{eqnarray*}

For $g_{ij}^0\wedge r_{jk}^1,\ w=L_{i(0)}H_{j(1)}H_k$, we have
\begin{eqnarray*}
&&(g_{ij}^0, r_{jk}^1)_w\\
&=&(L_{i(0)}H_j+H_{j(1)}DL_i-2H_{j(0)}L_i
-DH_{i+j})_{(1)}H_k-L_{i(0)}(H_{j(1)}H_k)\\
&=&-H_{j(1)}(L_{i(0)}H_k)+H_{i+j(0)}H_k\\
&\equiv &H_{j(1)}(H_{k(1)}DL_i-2H_{k(0)}L_i-DH_{i+k})+H_{i+j(0)}H_k\\
&\equiv &(H_{j(0)}H_{k})_{(2)}DL_i-2(H_{j(0)}H_{k})_{(1)}L_i
-H_{j(0)}H_{i+k}+H_{i+j(0)}H_k\\
&\equiv &-H_{j(0)}H_{i+k}+H_{i+j(0)}H_k\\
&\equiv &0\ mod(S_{_1}).
\end{eqnarray*}

For $g_{ij}^1\wedge q_{jkm}^0,\ w=L_{i(1)}H_{j(0)}L_{k+m}$, we have
\begin{eqnarray*}
&&(g_{ij}^1, q_{jkm}^0)_w\\
&=&(L_{i(1)}H_{j}+H_{j(1)}L_i-H_{i+j})_{(0)}L_{k+m}
-L_{i(1)}(H_{j(0)}L_{k+m}+H_{k(0)}L_{j+m}\\
&&-H_{0(0)}L_{j+k+m}-H_{j+k(0)}L_{m})\\
&=&-(L_{i(0)}H_{j})_{(1)}L_{k+m}+(H_{j(1)}L_i)_{(0)}L_{k+m}-H_{i+j(0)}L_{k+m}
-(L_{i(1)}H_{k})_{(0)}L_{j+m}\\
&&-(L_{i(0)}H_{k})_{(1)}L_{j+m}+(L_{i(1)}H_{0})_{(0)}L_{j+k+m}
+(L_{i(0)}H_{0})_{(1)}L_{j+k+m}\\
&&+(L_{i(1)}H_{j+k})_{(0)}L_{m}+(L_{i(0)}H_{j+k})_{(1)}L_{m}\\
&\equiv &-H_{j(0)}L_{i+k+m}-H_{k(0)}L_{i+j+m}+H_{0(0)}L_{i+j+k+m}+H_{j+k(0)}L_{i+m}\\
&\equiv &0\ mod(S_{_1}).
\end{eqnarray*}

For $g_{i,j+k}^1\wedge q_{jkm}^0,\ w=L_{i(1)}H_{j+k(0)}L_m$, we have
\begin{eqnarray*}
&&(g_{i,j+k}^1, -q_{jkm}^0)_w\\
&=&(L_{i(1)}H_{j+k}+H_{j+k(1)}L_i-H_{i+j+k})_{(0)}L_m
-L_{i(1)}(H_{j+k(0)}L_{m}-H_{j(0)}L_{k+m}\\
&&-H_{k(0)}L_{j+m}+H_{0(0)}L_{j+k+m})\\
&=&-(L_{i(0)}H_{j+k})_{(1)}L_{m}+(H_{j+k(1)}L_i)_{(0)}L_{m}-H_{i+j+k(0)}L_{m}
+(L_{i(1)}H_{k})_{(0)}L_{j+m}\\
&&+(L_{i(0)}H_{k})_{(1)}L_{j+m}-(L_{i(1)}H_{0})_{(0)}L_{j+k+m}
-(L_{i(0)}H_{0})_{(1)}L_{j+k+m}\\
&&+(L_{i(1)}H_{j})_{(0)}L_{k+m}+(L_{i(0)}H_{j})_{(1)}L_{k+m}\\
&\equiv &-H_{j+k(0)}L_{i+m}+H_{k(0)}L_{i+j+m}-H_{0(0)}L_{i+j+k+m}+H_{j(0)}L_{i+k+m}\\
&\equiv &0\ mod(S_{_1}).
\end{eqnarray*}

For $g_{ij}^1\wedge q_{jk}^1,\ w=L_{i(1)}H_{j(1)}L_k$, we have
\begin{eqnarray*}
&&(g_{ij}^1, q_{jk}^1)_w\\
&=&(L_{i(1)}H_j+H_{j(1)}L_i-H_{i+j})_{(1)}L_k
-L_{i(1)}(H_{j(1)}L_k-H_{0(1)}L_{j+k})\\
&=&H_{j(1)}(L_{i(1)}L_k)-H_{i+j(1)}L_k+(L_{i(1)}H_{0})_{(1)}L_{j+k}\\
&\equiv &H_{j(1)}L_{i+k}-H_{i+j(1)}L_k-(H_{0(1)}L_{i})_{(1)}L_{j+k}
+H_{i(1)}L_{j+k}\\
&\equiv &H_{j(1)}L_{i+k}-H_{0(1)}L_{i+j+k}-H_{i+j(1)}L_k
+H_{i(1)}L_{j+k}\\
&\equiv &0\ mod(S_{_1}).
\end{eqnarray*}

For $g_{ij}^1\wedge r_{jk}^0,\ w=L_{i(1)}H_{j(0)}H_k$, we have
\begin{eqnarray*}
&&(g_{ij}^1, r_{jk}^0)_w\\
&=&(L_{i(1)}H_j+H_{j(1)}L_i-H_{i+j})_{(0)}H_k-
L_{i(1)}(H_{j(0)}H_k-H_{0(0)}H_{j+k})\\
&=&-L_{i(0)}(H_{j(1)}H_k)+H_{j(1)}(L_{i(0)}H_k)-H_{j(0)}(L_{i(1)}H_k)
-H_{i+j(0)}H_k\\
&&+(L_{i(1)}H_{0})_{(0)}H_{j+k}+(L_{i(0)}H_{0})_{(1)}H_{j+k}\\
&\equiv &-H_{j(1)}(H_{k(1)}DL_i-2H_{k(0)}L_i-DH_{i+k})
+H_{j(0)}(H_{k(1)}L_i-H_{i+k})\\
&&-H_{i+j(0)}H_k-(H_{0(1)}L_{i})_{(0)}H_{j+k}+H_{i(0)}H_{j+k}\\
&\equiv &-(H_{j(0)}H_{k})_{(2)}DL_i+3(H_{j(0)}H_{k})_{(1)}L_i
+H_{j(1)}DH_{i+k}-H_{j(0)}H_{i+k}\\
&&-H_{i+j(0)}H_k-H_{0(1)}(L_{i(0)}H_{j+k})
+H_{0(0)}(L_{i(1)}H_{j+k})+H_{i(0)}H_{j+k}\\
&\equiv &(H_{0(0)}H_{j+k})_{(1)}L_i
+H_{0(1)}(H_{j+k(1)}DL_i-2H_{j+k(0)}L_i-DH_{i+j+k})\\
&&-H_{0(0)}(H_{j+k(1)}L_i-H_{i+j+k})\\
&\equiv &0\ mod(S_{_1}).
\end{eqnarray*}

For $g_{ij}^1\wedge r_{jk}^1,\ w=L_{i(1)}H_{j(1)}H_k$, we have
\begin{eqnarray*}
&&(g_{ij}^1, r_{jk}^1)_w\\
&=&(L_{i(1)}H_j+H_{j(1)}L_i-H_{i+j})_{(1)}H_k-
L_{i(1)}(H_{j(1)}H_k)\\
&=&H_{j(1)}(L_{i(1)}H_k)-H_{i+j(1)}H_k\\
&\equiv &-H_{j(1)}(H_{k(1)}L_i-H_{i+k})\\
&\equiv &0\ mod(S_{_1}).
\end{eqnarray*}

For $q_{ijk}^0\wedge s_{j+k,m}^0,\ w=H_{i(0)}L_{j+k(0)}L_m$, we have
\begin{eqnarray*}
&&(q_{ijk}^0, s_{j+k,m}^0)_w\\
&=&(H_{i(0)}L_{j+k}+H_{j(0)}L_{i+k}-H_{i+j(0)}L_{k}-H_{0(0)}L_{i+j+k})_{(0)}L_m\\
&&-H_{i(0)}(L_{j+k(0)}L_m-L_{0(0)}L_{j+k+m})\\
&=&H_{j(0)}(L_{i+k(0)}L_m)-H_{i+j(0)}(L_{k(0)}L_m)-H_{0(0)}(L_{i+j+k(0)}L_m)\\
&&+H_{i(0)}(L_{0(0)}L_{j+k+m})\\
&\equiv &(H_{j(0)}L_{i-j}-H_{i+j(0)}L_{-j}-H_{0(0)}L_{i}+H_{i(0)}L_{0})_{(0)}L_{j+k+m}\\
&\equiv &0\ mod(S_{_1}).
\end{eqnarray*}

For $q_{ijk}^0\wedge s_{km}^0,\ w=H_{i+j(0)}L_{k(0)}L_m$, we have
\begin{eqnarray*}
&&(-q_{ijk}^0, s_{km}^0)_w\\
&=&(H_{i+j(0)}L_{k}-H_{i(0)}L_{j+k}-H_{j(0)}L_{i+k}+H_{0(0)}L_{i+j+k})_{(0)}L_m\\
&&-H_{i+j(0)}(L_{k(0)}L_m-L_{0(0)}L_{k+m})\\
&=&-H_{j(0)}(L_{i+k(0)}L_m)-H_{i(0)}(L_{j+k(0)}L_m)+H_{0(0)}(L_{i+j+k(0)}L_m)\\
&&+H_{i+j(0)}(L_{0(0)}L_{k+m})\\
&\equiv &(-H_{j(0)}L_{i}-H_{i(0)}L_{j}+H_{0(0)}L_{i+j}+H_{i+j(0)}L_{0})_{(0)}L_{k+m}\\
&\equiv &0\ mod(S_{_1}).
\end{eqnarray*}

For $q_{ijk}^0\wedge s_{j+k,m}^1,\ w=H_{i(0)}L_{j+k(1)}L_m$, we have
\begin{eqnarray*}
&&(q_{ijk}^0, s_{j+k,m}^1)_w\\
&=&(H_{i(0)}L_{j+k}+H_{j(0)}L_{i+k}-H_{i+j(0)}L_{k}-H_{0(0)}L_{i+j+k})_{(1)}L_m\\
&&-H_{i(0)}(L_{j+k(1)}L_m-L_{j+k+m})\\
&\equiv &H_{j(0)}L_{i+k+m}-H_{i+j(0)}L_{k+m}-H_{0(0)}L_{i+j+k+m}+H_{i(0)}L_{j+k+m}\\
&\equiv &0\ mod(S_{_1}).
\end{eqnarray*}

For $q_{ijk}^0\wedge s_{km}^1,\ w=H_{i+j(0)}L_{k(1)}L_m$, we have
\begin{eqnarray*}
&&(-q_{ijk}^0, s_{km}^1)_w\\
&=&(H_{i+j(0)}L_{k}-H_{i(0)}L_{j+k}-H_{j(0)}L_{i+k}+H_{0(0)}L_{i+j+k})_{(1)}L_m\\
&&-H_{i+j(0)}(L_{k(1)}L_m-L_{k+m})\\
&\equiv &-H_{j(0)}L_{i+k+m}-H_{i(0)}L_{j+k+m}+H_{0(0)}L_{i+j+k+m}+H_{i+j(0)}L_{k+m}\\
&\equiv &0\ mod(S_{_1}).
\end{eqnarray*}

For $q_{ijk}^0\wedge g_{j+k,m}^0,\ w=H_{i(0)}L_{j+k(0)}H_m$, we have
\begin{eqnarray*}
&&(q_{ijk}^0, g_{j+k,m}^0)_w\\
&=&(H_{i(0)}L_{j+k}+H_{j(0)}L_{i+k}-H_{i+j(0)}L_{k}-H_{0(0)}L_{i+j+k})_{(0)}H_m\\
&&-H_{i(0)}(L_{j+k(0)}H_m+H_{m(1)}DL_{j+k}-2H_{m(0)}L_{j+k}-DH_{j+k+m})\\
&=&H_{j(0)}(L_{i+k(0)}H_m)-H_{i+j(0)}(L_{k(0)}H_m)-H_{0(0)}(L_{i+j+k(0)}H_m)\\
&&-H_{i(0)}(H_{m(1)}DL_{j+k})+2H_{i(0)}(H_{m(0)}L_{j+k})+H_{i(0)}DH_{j+k+m}\\
&\equiv &H_{j(0)}(H_{m(0)}L_{i+k})-H_{i+j(0)}(H_{m(0)}L_{k})
-H_{0(0)}(H_{m(0)}L_{i+j+k})\\
&&+H_{i(0)}(H_{m(0)}L_{j+k})\\
&\equiv &H_{0(0)}(H_{i+m(0)}L_{j+k}+H_{j+m(0)}L_{i+k}-H_{i+j+m(0)}L_{k}\\
&&-H_{m(0)}L_{i+j+k})\\
&\equiv & \left\{
\begin{aligned}
&H_{0(0)}(H_{0(0)}L_{j+k+i+m}+H_{j+m(0)}L_{i+k}-H_{j(0)}L_{k+i+m}\\
&\ \ \ \ \ \ \ \ \ \ \ \ \ \   -H_{m(0)}L_{i+j+k}),\ \ \ \ \ \  i+m=0   \\
&H_{0(0)}(-H_{j(0)}L_{k+i+m}+H_{0(0)}L_{i+j+k+m}+H_{j+m(0)}L_{k+i}\\
&\ \ \ \ \ \ \ \ \ \ \ \ \ \  -H_{m(0)}L_{i+j+k}),\ \ \ \ \ \   i+m\neq 0
\end{aligned}
\right. \\
&\equiv &0\ mod(S_{_1}).
\end{eqnarray*}

For $q_{ijk}^0\wedge g_{km}^0,\ w=H_{i+j(0)}L_{k(0)}H_m$, we have
\begin{eqnarray*}
&&(-q_{ijk}^0, g_{km}^0)_w\\
&=&(H_{i+j(0)}L_{k}-H_{i(0)}L_{j+k}-H_{j(0)}L_{i+k}+H_{0(0)}L_{i+j+k})_{(0)}H_m\\
&&-H_{i+j(0)}(L_{k(0)}H_m+H_{m(1)}DL_{k}-2H_{m(0)}L_{k}-DH_{k+m})\\
&=&-H_{j(0)}(L_{i+k(0)}H_m)-H_{i(0)}(L_{j+k(0)}H_m)+H_{0(0)}(L_{i+j+k(0)}H_m)\\
&&-H_{i+j(0)}(H_{m(1)}DL_{k})+2H_{i+j(0)}(H_{m(0)}L_{k})+H_{i+j(0)}DH_{k+m}\\
&\equiv &-H_{j(0)}(H_{m(0)}L_{i+k})-H_{i(0)}(H_{m(0)}L_{j+k})
+H_{0(0)}(H_{m(0)}L_{i+j+k})\\
&&+H_{i+j(0)}(H_{m(0)}L_{k})\\
&\equiv &0\ mod(S_{_1}).
\end{eqnarray*}

For $q_{ijk}^0\wedge g_{j+k,m}^1,\ w=H_{i(0)}L_{j+k(1)}H_m$, we have
\begin{eqnarray*}
&&(q_{ijk}^0, g_{j+k,m}^1)_w\\
&=&(H_{i(0)}L_{j+k}+H_{j(0)}L_{i+k}-H_{i+j(0)}L_{k}-H_{0(0)}L_{i+j+k})_{(1)}H_m\\
&&-H_{i(0)}(L_{j+k(1)}H_m+H_{m(1)}L_{j+k}-H_{j+k+m})\\
&=&H_{j(0)}(L_{i+k(1)}H_m)-H_{i+j(0)}(L_{k(1)}H_m)-H_{0(0)}(L_{i+j+k(1)}H_m)\\
&&-H_{i(0)}(H_{m(1)}L_{j+k})+H_{i(0)}H_{j+k+m}\\
&\equiv &-H_{j(0)}(H_{m(1)}L_{i+k})+H_{j(0)}H_{i+k+m}+H_{i+j(0)}(H_{m(1)}L_{k})+H_{i+j(0)}H_{m+k}\\
&&+H_{0(0)}(H_{m(1)}L_{i+j+k})+H_{0(0)}H_{i+j+k+m}-H_{i(0)}(H_{m(1)}L_{j+k})+H_{i(0)}H_{j+k+m}\\
&\equiv &0\ mod(S_{_1}).
\end{eqnarray*}

For $q_{ijk}^0\wedge g_{km}^1,\ w=H_{i+j(0)}L_{k(1)}H_m$, we have
\begin{eqnarray*}
&&(-q_{ijk}^0, g_{km}^1)_w\\
&=&(H_{i+j(0)}L_{k}-H_{i(0)}L_{j+k}-H_{j(0)}L_{i+k}+H_{0(0)}L_{i+j+k})_{(1)}H_m\\
&&-H_{i+j(0)}(L_{k(1)}H_m+H_{m(1)}L_k-H_{k+m})\\
&=&-H_{i(0)}(L_{j+k(1)}H_m)-H_{j(0)}(L_{i+k(1)}H_m)+H_{0(0)}(L_{i+j+k(1)}H_m)\\
&&-H_{i+j(0)}(H_{m(1)}L_{k})+H_{i+j(0)}H_{k+m}\\
&\equiv &0\ mod(S_{_1}).
\end{eqnarray*}

For $q_{ij}^1\wedge s_{jk}^0,\ w=H_{i(1)}L_{j(0)}L_k$, we have
\begin{eqnarray*}
&&(q_{ij}^1, s_{jk}^0)_w\\
&=&(H_{i(1)}L_j-H_{0(1)}L_{i+j})_{(0)}L_k-H_{i(1)}(L_{j(0)}L_k-L_{0(0)}L_{j+k})\\
&=&-H_{i(0)}(L_{j(1)}L_k)-H_{0(1)}(L_{i+j(0)}L_k)
+H_{0(0)}(L_{i+j(1)}L_k)+(H_{i(1)}L_{0})_{(0)}L_{j+k}\\
&&+H_{i(0)}(L_{0(1)}L_{j+k})\\
&\equiv &-H_{i(0)}L_{j+k}-H_{0(1)}(L_{0(0)}L_{i+j+k})+H_{0(0)}L_{i+j+k}
+(H_{0(1)}L_{i})_{(0)}L_{j+k}\\
&&+H_{i(0)}L_{j+k}\\
&\equiv &-H_{0(1)}(L_{0(0)}L_{i+j+k})+H_{0(0)}L_{i+j+k}
+H_{0(1)}(L_{i(0)}L_{j+k})-H_{0(0)}L_{i+j+k}\\
&\equiv &0\ mod(S_{_1}).
\end{eqnarray*}

For $q_{ij}^1\wedge s_{jk}^1,\ w=H_{i(1)}L_{j(1)}L_k$, we have
\begin{eqnarray*}
&&(q_{ij}^1, s_{jk}^1)_w\\
&=&(H_{i(1)}L_j-H_{0(1)}L_{i+j})_{(1)}L_k
-H_{i(1)}(L_{j(1)}L_k-L_{j+k})\\
&=&-H_{0(1)}(L_{i+j(1)}L_k)+H_{i(1)}L_{j+k}\\
&\equiv &-H_{0(1)}L_{i+j+k}+H_{i(1)}L_{j+k}\\
&\equiv &0\ mod(S_{_1}).
\end{eqnarray*}

For $q_{ij}^1\wedge g_{jk}^0,\ w=H_{i(1)}L_{j(0)}H_k$, we have
\begin{eqnarray*}
&&(q_{ij}^1, g_{jk}^0)_w\\
&=&(H_{i(1)}L_j-H_{0(1)}L_{i+j})_{(0)}H_k-H_{i(1)}
(L_{j(0)}H_k+H_{k(1)}DL_j-2H_{k(0)}L_j-DH_{j+k})\\
&=&-H_{i(0)}(L_{j(1)}H_k)-H_{0(1)}(L_{i+j(0)}H_k)+H_{0(0)}
(L_{i+j(1)}H_k)-H_{i(1)}(H_{k(1)}DL_j) \\
&&+2H_{i(1)}(H_{k(0)}L_j)+H_{i(1)}DH_{j+k} \\
&\equiv &-H_{i(0)}(L_{j(1)}H_k)-H_{0(1)}(L_{i+j(0)}H_k)
+H_{0(0)}(L_{i+j(1)}H_k)-H_{i(0)}(H_{k(2)}DL_j) \\
&&+2H_{i(0)}(H_{k(1)}L_j)+H_{i(0)}H_{j+k} \\
&\equiv &H_{i(0)}(H_{k(1)}L_j)-H_{i(0)}H_{j+k}
+H_{0(1)}(H_{k(1)}DL_{i+j}-2H_{k(0)}L_{i+j}-DH_{i+j+k})\\
&&-H_{0(0)}(H_{k(1)}L_{i+j})
+H_{0(0)}H_{i+j+k}+H_{i(0)}H_{j+k} \\
&\equiv &H_{0(1)}(H_{k(1)}DL_{i+j})
-2H_{0(1)}(H_{k(0)}L_{i+j})\\
&\equiv &H_{0(0)}(H_{k(2)}DL_{i+j})-2H_{0(0)}(H_{k(1)}L_{i+j})\\
&\equiv &0\ mod(S_{_1}).
\end{eqnarray*}

For $q_{ij}^1\wedge g_{jk}^1,\ w=H_{i(1)}L_{j(1)}H_k$, we have
\begin{eqnarray*}
&&(q_{ij}^1, g_{jk}^1)_w\\
&=&(H_{i(1)}L_{j}-H_{0(1)}L_{i+j})_{(1)}H_k
-H_{i(1)}(L_{j(1)}H_k+H_{k(1)}L_j-H_{j+k})\\
&=&-H_{0(1)}(L_{i+j(1)}H_k)-H_{i(1)}(H_{k(1)}L_j)+H_{i(1)}H_{j+k}\\
&\equiv &H_{0(1)}(H_{k(1)}L_{i+j})-H_{0(1)}H_{i+j+k}\\
&\equiv &0\ mod(S_{_1}).
\end{eqnarray*}

For $r_{ij}^0\wedge q_{jkm}^0,\ w=H_{i(0)}H_{j(0)}L_{k+m}$, we have
\begin{eqnarray*}
&&(r_{ij}^0, q_{jkm}^0)_w\\
&=&(H_{i(0)}H_j-H_{0(0)}H_{i+j})_{(0)}L_{k+m}-H_{i(0)}(H_{j(0)}L_{k+m}+H_{k(0)}L_{j+m}\\
&&-H_{j+k(0)}L_{m}-H_{0(0)}L_{j+k+m})\\
&\equiv &-H_{0(0)}(H_{i+j(0)}L_{k+m})-H_{0(0)}(H_{i+k(0)}L_{j+m})+H_{0(0)}(H_{i+j+k(0)}L_{m})\\
&&+H_{0(0)}(H_{i(0)}L_{j+k+m})\\
&\equiv &0\ mod(S_{_1}).
\end{eqnarray*}

For $r_{i,j+k}^0\wedge q_{jkm}^0,\ w=H_{i(0)}H_{j+k(0)}L_{m}$, we have
\begin{eqnarray*}
&&(r_{i,j+k}^0, -q_{jkm}^0)_w\\
&=&(H_{i(0)}H_{j+k}-H_{0(0)}H_{i+j+k})_{(0)}L_{m}-H_{i(0)}(H_{j+k(0)}L_{m}-H_{j(0)}L_{k+m}\\
&&-H_{k(0)}L_{j+m}+H_{0(0)}L_{j+k+m})\\
&\equiv &-H_{0(0)}(H_{i+j+k(0)}L_{m})+H_{0(0)}(H_{i+j(0)}L_{k+m})+H_{0(0)}(H_{i+k(0)}L_{j+m})\\
&&-H_{0(0)}(H_{i(0)}L_{j+k+m})\\
&\equiv &0\ mod(S_{_1}).
\end{eqnarray*}

For $r_{ij}^0\wedge q_{jk}^1,\ w=H_{i(0)}H_{j(1)}L_k$, we have
\begin{eqnarray*}
&&(r_{ij}^0, q_{jk}^1)_w\\
&=&(H_{i(0)}H_j-H_{0(0)}H_{i+j})_{(1)}L_k-H_{i(0)}(H_{j(1)}L_k-H_{0(1)}L_{j+k})\\
&=&-H_{0(0)}(H_{i+j(1)}L_k)+(H_{i(0)}H_{0})_{(1)}L_{j+k}\\
&\equiv &-H_{0(0)}(H_{i+j(1)}L_k)+H_{0(0)}(H_{i(1)}L_{j+k})\\
&\equiv &0\ mod(S_{_1}).
\end{eqnarray*}

For $r_{ij}^1\wedge q_{jkm}^0,\ w=H_{i(1)}H_{j(0)}L_{k+m}$, we have
\begin{eqnarray*}
&&(r_{ij}^1, q_{jkm}^0)_{w_1}\\
&=&(H_{i(1)}H_j)_{(0)}L_{k+m}-H_{i(1)}(H_{j(0)}L_{k+m}+H_{k(0)}L_{j+m}
-H_{j+k(0)}L_{m}\\
&&-H_{0(0)}L_{j+k+m})\\
&=&-H_{i(0)}(H_{j(1)}L_{k+m})-(H_{i(1)}H_k)_{(0)}L_{j+m}-H_{i(0)}(H_{k(1)}L_{j+m})\\
&&+(H_{i(1)}H_{j+k})_{(0)}L_{m}+H_{i(0)}(H_{j+k(1)}L_{m})+(H_{i(1)}H_0)_{(0)}L_{j+k+m}\\
&&+H_{i(0)}(H_{0(1)}L_{j+k+m})\\
&\equiv &0\ mod(S_{_1}).
\end{eqnarray*}

For $r_{i,j+k}^1\wedge q_{jkm}^0,\ w=H_{i(1)}H_{j+k(0)}L_{m}$, we have
\begin{eqnarray*}
&&(r_{i,j+k}^1, -q_{jkm}^0)_{w_2}\\
&=&(H_{i(1)}H_{j+k})_{(0)}L_{m}+H_{i(1)}(H_{j(0)}L_{k+m}+H_{k(0)}L_{j+m}
-H_{j+k(0)}L_{m}\\
&&-H_{0(0)}L_{j+k+m})\\
&=&-H_{i(0)}(H_{j+k(1)}L_{m})+(H_{i(1)}H_j)_{(0)}L_{k+m}
+H_{i(0)}(H_{j(1)}L_{k+m})\\
&&+(H_{i(1)}H_k)_{(0)}L_{j+m}+H_{i(0)}(H_{k(1)}L_{j+m})-(H_{i(1)}H_0)_{(0)}L_{j+k+m}\\
&&-H_{i(0)}(H_{0(1)}L_{j+k+m})\\
&\equiv &0\ mod(S_{_1}).
\end{eqnarray*}

For $r_{ij}^1\wedge q_{jk}^1,\ w=H_{i(1)}H_{j(1)}L_k$, we have
\begin{eqnarray*}
(r_{ij}^1, q_{jk}^1)_w
&=&(H_{i(1)}H_j)_{(1)}L_k-H_{i(1)}(H_{j(1)}L_k-H_{0(1)}L_{j+k})\\
&=&(H_{i(1)}H_{0})_{(1)}L_{j+k}\\
&\equiv &0\ mod(S_{_1}).
\end{eqnarray*}

For any $x\in\{L_i, H_j\}_{i,j\in\mathbb{Z}}, n\geq 2$, we have
\begin{eqnarray*}
x_{(n)}s_{ij}^0
&=&x_{(n)}(L_{i(0)}L_j-L_{0(0)}L_{i+j})\\
&=&-\sum_{t\geq 1}(-1)^{t}\binom{n}{t}x_{(n-t)}(L_{i(t)}L_j-L_{0(t)}L_{i+j})\\
&\equiv &nx_{(n-1)}(L_{i(1)}L_j-L_{0(1)}L_{i+j})\\
&\equiv& 0\ mod(s_{_1}),\\
x_{(n)}s_{ij}^1
&=&x_{(n)}(L_{i(1)}L_j-L_{i+j})\\
&=&-\sum_{t\geq 1}(-1)^{t}\binom{n}{t}x_{(n-t)}(L_{i(1+t)}L_j)\\
&\equiv &0\ mod(S_{_1}),\\
x_{(n)}r_{ij}^0
&=&x_{(n)}(H_{i(0)}H_j-H_{0(0)}H_{i+j})\\
&=&=-\sum_{t\geq 1}(-1)^{t}\binom{n}{t}x_{(n-t)}(H_{i(t)}H_j-H_{0(t)}H_{i+j})\\
&=&nx_{(n-1)}(H_{i(1)}H_j-H_{0(1)}H_{i+j})\\
&\equiv& 0\ mod(S_{_1}),\\
\end{eqnarray*}

\begin{eqnarray*}
x_{(n)}r_{ij}^1
&=&x_{(n)}(H_{i(1)}H_j)\\
&=&-\sum_{t\geq 1}(-1)^{t}\binom{n}{t}x_{(n-t)}(H_{i(1+t)}H_j)\\
&\equiv &0\ mod(S_{_1}),\\
x_{(n)}g_{ij}^1
&=&x_{(n)}(L_{i(1)}H_j+H_{j(1)}L_i-H_{i+j})\\
&=&-\sum_{t\geq 1}(-1)^{t}\binom{n}{t}x_{(n-t)}(L_{i(1+t)}H_j+H_{j(1+t)}L_i)\\
&\equiv &0\ mod(S_{_1}),\\
x_{(n)}q_{ijk}^0
&=&x_{(n)}(H_{i(0)}L_{j+k}+H_{j(0)}L_{i+k}
-H_{i+j(0)}L_{k}-H_{0(0)}L_{i+j+k})\\
&=&-\sum_{t\geq 1}(-1)^{t}\binom{n}{t}x_{(n-t)}(H_{i(t)}L_{j+k}+H_{j(t)}L_{i+k}
-H_{i+j(t)}L_{k}\\
&&-H_{0(t)}L_{i+j+k})\\
&=&nx_{(n-1)}(H_{i(1)}L_{j+k}+H_{j(1)}L_{i+k}
-H_{i+j(1)}L_{k}-H_{0(1)}L_{i+j+k})\\
&\equiv &0\ mod(S_{_1}),\\
x_{(n)}q_{ij}^1
&=&x_{(n)}(H_{i(1)}L_j-H_{0(1)}L_{i+j})\\
&=&-\sum_{t\geq 1}(-1)^{t}\binom{n}{t}x_{(n-t)}(H_{i(1+t)}L_j-H_{0(1+t)}L_{i+j})\\
&\equiv &0\ mod(S_{_1}),\\
x_{(n)}g_{ij}^0
&=&x_{(n)}(L_{i(0)}H_j+H_{j(1)}DL_i-2H_{j(0)}L_i-DH_{i+j})\\
&=&-\sum_{t\geq 1}(-1)^{t}\binom{n}{t}x_{(n-t)}(L_{i(t)}H_j+H_{j(1+t)}DL_i-2H_{j(t)}L_i)\\
&&-nx_{(n-1)}H_{i+j}\\
&=&nx_{(n-1)}(L_{i(1)}H_j+H_{j(2)}DL_i-2H_{j(1)}L_i)-nx_{(n-1)}H_{i+j} \\
&=& \left\{
\begin{aligned}
& 0,\ \  n>2   \\
&2x_{(1)}(L_{i(1)}H_j+H_{j(2)}DL_i-2H_{j(1)}L_i)-2x_{(1)}H_{i+j},\ \   n=2,
\end{aligned}
\right. \\
g_{ij(n)}^0x
&=&(L_{i(0)}H_j+H_{j(1)}DL_i-2H_{j(0)}L_i-DH_{i+j})_{(n)}x\\
&=&-nH_{j(1)}(L_{i(n-1)}x)+(n+1)H_{j(0)}(L_{i(n)}x)+nH_{i+j(n-1)}x\\
&=&-nH_{j(1)}(L_{i(n-1)}x)+nH_{i+j(n-1)}x\\
&=& \left\{
\begin{aligned}
& 0,\ \  n>2   \\
&-2H_{j(1)}(L_{i(1)}x)+2H_{i+j(1)}x,\ \   n=2.
\end{aligned}
\right.
\end{eqnarray*}

Then for any $k\in\mathbb{Z}$, we have
\begin{eqnarray*}
g_{ij(2)}^0L_k&\equiv& -2H_{j(1)}L_{i+k}+2H_{i+j(1)}L_k\\
&\equiv &0\ mod(S_{_1}),\\
g_{ij(2)}^0H_k&=& -2H_{j(1)}(L_{i(1)}H_k)+2H_{i+j(1)}H_k\\
&\equiv& 2H_{j(1)}(H_{k(1)}L_i)-2H_{j(1)}H_{i+k}\\
&\equiv&0\ mod(S_{_1}),\\
L_{k(2)}g_{ij}^0
&=&2L_{k(1)}(L_{i(1)}H_j+H_{j(2)}DL_i-2H_{j(1)}L_i)-2L_{k(1)}H_{i+j}\\
&=&2(L_{k(1)}L_{i})_{(1)}H_j+2(L_{k(1)}H_{j})_{(2)}DL_i-4(L_{k(1)}H_{j})_{(1)}L_i-2L_{k(1)}H_{i+j}\\
&\equiv&2L_{i+k(1)}H_j-2(H_{j(1)}L_{k})_{(2)}DL_i+2H_{j+k(2)}DL_i+4(H_{j(1)}L_{k})_{(1)}L_i\\
&&-4H_{j+k(1)}L_i+2H_{i+j(1)}L_k-2H_{i+j+k}\\
&\equiv&-2H_{j(1)}L_{i+k}+2H_{i+j(1)}L_k\\
&\equiv&0\ mod(S_{_1}),\\
H_{k(2)}g_{ij}^0&=&2H_{k(1)}(L_{i(1)}H_j+H_{j(2)}DL_i-2H_{j(1)}L_i)-2H_{k(1)}H_{i+j}\\
&\equiv&-2H_{k(1)}(H_{j(1)}L_i)+2H_{k(1)}H_{i+j}\\
&\equiv&0\ mod(S_{_1}).
\end{eqnarray*}

Therefore
 $S_{_1}$ is a Gr\"{o}bner-Shirshov basis in
 $C(\{L_i, H_j\}_{i,j\in \mathbb{Z}},\ N=2)$.

(iii)\ This part follows from Theorem \ref{t1} and  Proposition \ref{t2}.
\ \ \ \  $\square$

\ \

Not each Lie conformal algebra is  embeddable into its universal
enveloping associative conformal algebra, see, for example, \cite{Ro00}.

\ \

\noindent{\bf Remark}\ In the Examples \ref{exam1} and \ref{exam2}, if
we replace the complex field $\mathbb{C}$ with an arbitrary field of characteristic $0$, then all results hold.

\end{document}